\RequirePackage{fix-cm}
\documentclass[smallextended]{svjour3}       
\smartqed  
\usepackage{graphicx,comment}
\usepackage{amssymb}
\usepackage{amsmath,color,latexsym} 
\usepackage{multirow}
\usepackage[table,xcdraw]{xcolor}
\newcommand{\eps}{\varepsilon}

\usepackage[normalem]{ulem}

\begin{document}

\title{
Finite-time stability of 
polyhedral sweeping processes with application to elastoplastic systems
\thanks{The first and second authors were supported by the National Science Foundation grant CMMI-1916876.}
}

\titlerunning{Finite-time stability of 
polyhedral sweeping processes}        

\author{Ivan Gudoshnikov \and
        Oleg Makarenkov \and Dmitry Rachinskiy 
}


\institute{I. Gudoshnikov \at
Department of Physics, Arizona State University, Tempe AZ 85281\\
      \email{Ivan.Gudoshnikov@asu.edu}      
     \and O. Makarenkov, D. Rachinskiy \at  
Department of Mathematical Sciences, University of Texas at Dallas,
Richardson TX 75080       \\
              \email{makarenkov@utdallas.edu},
             \email{dmitry.rachinskiy@utdallas.edu} 
}

\date{Received: date / Accepted: date}

\maketitle

\begin{abstract} We use the ideas of Adly-Attoych-Cabot [Adv. Mech. Math., 12, Springer, 2006] on finite-time stabilization of dry friction oscillators to establish a theorem on finite-time stabilization of differential inclusions with a moving polyhedral constraint (known as polyhedral sweeping processes) of the form $C+c(t).$ We then employ the ideas of Moreau [New variational techniques in mathematical physics, CIME, 1973] to apply our theorem to a system of elastoplastic springs with a displacement-controlled loading. We show that verifying the condition of the theorem ultimately leads to the following two problems: (i) identifying the active vertex ``A'' or the active face ``A'' of the polyhedron that the vector $c'(t)$ points at; (ii) computing the distance from $c'(t)$ to the normal cone to the polyhedron at ``A''. We provide a computational guide to implement steps (i)-(ii) in the case of an arbitrary elastoplastic system and apply the guide to a particular example. Due to the simplicity of the particular example, we can solve (i)-(ii) by the methods of linear algebra and minor combinatorics. 

\vskip0.2cm

\end{abstract}

\section{Introduction} \label{intro}


Finite-time stability of an attractor is typical for differential equations with nonsmooth right-hand-sides. This fact is used in control theory since long ago. Finite-time stability in differential equations with nonsmooth right-hand-sides is often proved by showing that a Lyapunov function $V$ satisfies the estimate (see e.g. 
Bernuau et al \cite{Bernuau-et-al}, Bhat-Bernstein \cite{bhat1}, Oza et al \cite{Oza-et-al}, Sanchez et al \cite{Sanchez-et-al})
\begin{equation}\label{lyap2}
   \frac{d}{dt}[V(x(t))]+2\varepsilon\sqrt{V(x(t))}\le 0,\quad \mbox{a.e. on }[0,\infty),
\end{equation}
for some $\varepsilon>0$, where $x$ is a solution. Specifically, if (\ref{lyap2}) holds for a function $x(t)$, then $V(x(t_1))=0$ at some $0\le t_1,$ where (see Lemma~\ref{techlemma0}) 
\begin{equation}\label{t1}
t_1\le \dfrac{1}{\eps}V(x(0)).
\end{equation}

\vskip0.2cm

\noindent Motivated by applications in frictional mechanics, Adly et al \cite{Adly} extended the Lyapunov function approach to finite-time stability analysis of differential inclusions. Let $\nabla f(x)$ be the gradient of a function $f:\mathbb{R}\to\mathbb{R}^n$, $\partial \Phi(x)$ be the subdifferential of a convex function $\Phi:\mathbb{R}\to\mathbb{R}^n$, and $B_\varepsilon(0)$ be the ball of $\mathbb{R}^n$ of radius $\varepsilon$ centered at 0. By focusing on differential inclusions of the form 
\begin{equation}\label{Adlyinclusion}
   -\ddot x(t)-\nabla f(x(t))\in \partial\Phi(\dot x(t)), 
\end{equation}
the paper \cite{Adly} discovered (see the proof of \cite[Theorem~24.8]{Adly}) that the property
\begin{equation}\label{Beps}
   -\nabla f(x(t))+B_{\varepsilon}(0)\subset \partial\Phi(0), \quad \mbox{a.e. on }[0,\infty),
\end{equation}
implies (\ref{lyap2})
for a suitable Lyapunov function $V$ that measures the distance from $\dot x(t)$ to $0$ and for any solution $x$ of (\ref{Adlyinclusion}). 

\vskip0.2cm

\noindent  More recently, a significant interest in the study of finite-time stability of differential inclusions has been due to new applications in elastoplasticity (see e.g. Gudoshnikov et al \cite{GKMV}). We remind the reader that according to the pioneering work by Moreau \cite{Moreau} (see also Gudoshnikov-Makarenkov \cite{ESAIM}), the stresses in a network of $m$ elastoplastic springs with time-varying displacement-controlled loadings are governed by 
\begin{equation}\label{sp}
   -\dot y\in N^A_{C(t)}(y),\quad y\in \mathcal{V},\quad \begin{array}{l}\mathcal{V}\mbox{ is a }d-\mbox{dimensional subspace of }\mathbb{R}^m\\
  \mbox{with the scalar product}\ (x,y)_A=\left<x,Ay\right>,\end{array} 
\end{equation}
where $A$ is a positive diagonal $m\times m$-matrix, and $N_{C(t)}^A(y)$ is a normal cone to the set 
\begin{equation}
  C(t)=C+c(t),\quad  C=\bigcap_{j=1}^m \mathcal{V}_j,  \quad  
  \def\arraystretch{1.5}
  \begin{array}{rcl} \mathcal{V}_j&=&L(-1,j)\cap L(+1,j),\\ 
  L(\alpha,j)&=&\left\{y\in \mathcal{V}\hskip-0.1cm:\left<\alpha e_j,Ay\right>\le \alpha c_{j}^{\alpha}\right\},
  \end{array}  \label{C(t)}
\end{equation}
at a point $y$, with appropriate $d,c_j^-,$ $c_j^+,$ $c(t)$ that define mechanical parameters of the network of elastoplastic springs and the displacement-controlled loadings (see Section~\ref{mechanics}). The solutions $y(t)$ of differential inclusion (\ref{sp}) never escape from $C(t)$ (i.e. $y(t)$ is swept by $C(t)$)  for which reason (\ref{sp}) is called {\it sweeping process.} Spring $j$ undergoes {\it plastic deformation} when the inequality $c_j^-<\left<e_j,Ay\right><c_j^+$ is violated.  Therefore, knowledge of the evolution of $y(t)$ allows to make conclusions about the regions of plastic deformation (that lead to {\it low-cycle fatigue} or {\it incremental failure}, see Yu \cite[\S4.6]{19}). 

\vskip0.2cm

\noindent Krejci \cite{Krejci} proved that if $c(t)$ is $T$-periodic then the set $Y$ of $T$-periodic solutions of (\ref{sp}) is always asymptotically stable. Examining finite-time stability of $Y$ is a hard problem because it requires the knowledge of particular elements of $Y$ (our work in progress). Accordingly, conditions for finite-time stability of $Y$ are not going to be easily verifiable except for the case where $Y$ consists of just one solution (see Gudoshnikov et al \cite{GKMV}). 

\vskip0.2cm

\noindent At the same time, predicting the behavior of solutions of sweeping process (\ref{sp}) within a guaranteed time is of crucial importance for materials science. 
Current methods of computing the asymptotic response of networks of elastoplastic springs (see e.g. Boudy et al \cite{bouby}, Zouain-SantAnna \cite{brazil}) run the numeric routine until the difference between the responses corresponding to two successive cycles of  loading get smaller than a prescribed tolerance (without any estimate as for how soon such a desired accuracy will be reached).  

\vskip0.2cm

\noindent The present paper adapts condition (\ref{Beps}) in order to predict the behavior of solutions of (\ref{sp}) within a guaranteed finite time. We don't prove the  finite-time stability of $Y$, but still prove that all solutions of (\ref{sp}) will be confined within a certain  computable set.
 Specifically, let $F(t)$ be a facet of $C(t)$ and let ${\rm ri}(F)$ denote the relative interior of $F$. 
 This means that $F(t)$ can be expressed as  
\begin{equation}\label{F(t)} \def\arraystretch{1.5}
\hskip-0.1cm \begin{array}{l}
   F(t)=F+c(t), \\ 
   \displaystyle
   F=\left(\bigcap\limits_{(\alpha,j)\in I_0} \overline{\overline{L}}(\alpha,j)\right)\cap \left(\bigcap_{i=1}^M\bigcap\limits_{(\alpha,j)\in I_i} {{L}}(\alpha,j)\right),   \\
\overline{\overline{L}}(\alpha,j)=\left\{y\in \mathcal{V}\hskip-0.1cm:\left< e_j,Ay\right>=  c_{j}^{\alpha}\right\}, 
 \end{array}
\end{equation}
where the ingredients in (\ref{F(t)}) satisfy the following assumptions:
\begin{eqnarray}
&& \hskip-0.7cm \begin{array}{l}\mbox{each }I_0\cup I_i\mbox{ defines}\\
 \mbox{a vertex of }F:
 \end{array}\ \   \begin{array}{l}
\mbox{for each }i\in\overline{1,M},\ \{y:y\in\overline{\overline{L}}(\alpha,j),\ (\alpha,j)\in I_0\cup I_i\}\\
\mbox{is a non-empty singleton }\{y_{*,i}\},\ 
\end{array} \label{prop1}  \\
&&
\hskip-0.7cm\begin{array}{l}F\mbox{ has no other vertices,}\\
\mbox{all constraints are counted}:\end{array}\ \  
\begin{array}{l} \mbox{each }y\in F\mbox{ defines }J_y\subset\overline{1,M}\cup\emptyset,\mbox{ such}\\ \mbox{that }\{(\alpha,j)\hskip-0.06cm:y\in\overline{\overline L}(\alpha,j)\}=I_0\cup\bigcap_{i\in J_y} I_i,
\end{array}\label{noother}\\
&&\hskip-0.7cm \mbox{vertices do not coincide}:\ \ J_{y_{*,i}}=\{i\},\ i\in\overline{1,M},\label{notcoincide}\\
&&\hskip-0.7cm \mbox{vertices do not reduce ri}(F):\ \ \{(\alpha,j)\hskip-0.06cm:y\in\overline{\overline L}(\alpha,j)\}=I_0,\ y\in{\rm ri}(F),\label{notreduced}\\
&&\hskip-0.7cm F\mbox{ is feasible:} \ \  F\subset C, \label{Ffeasible} \\
&& \hskip-0.7cm \begin{array}{l}\overline{\overline{L}}(\alpha,j),\ (\alpha,j)\in I_0\cup I_i,\mbox{ are  independent}:
 \end{array}\ \   \begin{array}{l}
|I_0|+|I_i|=d,\ i\in\overline{1,M}.\end{array} \label{independent}
\end{eqnarray}
The case where $F$ is a just a vertex of $C$ is accounted for by $M=0,$ see formulas (\ref{Fred1})-(\ref{Ffeasible_}) below for the corresponding reduced form of (\ref{F(t)})-(\ref{independent}). 
We prove that if 
\begin{equation}\label{A0}
   -\dot c(t)+B_\eps^A(0)\cap N_F^A(y)\subset N_{C}^A(y),\quad y\in F, \ \mbox{a.a.}\ t\in[0,\tau_d],
\end{equation}
where $B_\eps^A(0)$ is a ball in the norm induced by the scalar product (\ref{sp}),
then, for any solution $y(t)$ of (\ref{sp}), the function
\begin{equation}\label{x(t)}
   x(t)=y(t)-c(t)
\end{equation}
 satisfies the estimate (\ref{lyap2}) on $[0,\tau_d]$ for a suitable Lyapunov function $V$ that measures the distance from $x(t)$ to $F.$ Since, by (\ref{F(t)}), the distance from $x(t)$ to $F$ equals the distance from $y(t)$ to $F(t)$, then the relation $\tau_d\ge\dfrac{1}{\eps}V(x(0))$ ensures one-period reachability of the facet $F(t)$ when $c(t)$ is $T$-periodic.

\vskip0.2cm

\noindent The paper is organized as follows. In Section~\ref{sec2} we prove our main result (Theorem~\ref{mainthm}) which uses condition (\ref{A0}) in order to estimate the time it takes for all solution of (\ref{sp}) to reach the facet $F(t)$. The proof of Theorem~\ref{mainthm} relies on the ideas of Adly et al \cite{Adly}, which are used in \cite{Adly} to establish finite-time stability of a frictional system. However, for the proof of Theorem~\ref{mainthm}, we reformulated the ideas of \cite{Adly} in terms of a suitable Lyapunov function, which can be of independent interest for applied sciences. The two successive sections (Section~\ref{sec3} and Section~\ref{sec4}) derive the corollaries of Theorem~\ref{mainthm} for the case where $F(t)$ consists of just one point (Corollary~\ref{cor1}) and where $F(t)$ is an entire facet (Corollary~\ref{thm3}). Furthermore, Section~\ref{sec4} discovers (Corollary~\ref{corforstep3}) that when the validity of condition (\ref{A0}) is known in the interior points of $F$ only, the convergence of all solutions of (\ref{sp}) to $F(t)$ still occurs in finite time, but the estimate of the time of convergence is no longer available. In order to establish Corollaries~\ref{thm3} and \ref{corforstep3}, Section~\ref{sec4} derives several computations formulas for normal cones to the moving constraint $C(t)$ and its facet $F(t)$ (formulas (\ref{goodformula2})-(\ref{goodformula1})), which take roots in Rockafellar-Wets \cite{Rockafellar-Wets} and which can be of independent interest for set-valued analysis. 

\vskip0.2cm

\noindent Section~\ref{mechanics} summarizes the Moreau approach \cite{Moreau} towards the use of sweeping process (\ref{sp}) to model networks of elastoplastic springs. The notations of this section follow Gudoshnikov-Makarenkov \cite{PhysicaD,ESAIM}. Section~\ref{guide} provides a step-by-step guide for application of the results of Sections~\ref{sec2}-\ref{sec4} to sweeping processes coming from networks of elastoplastic springs. In particular, Section~\ref{guide} uses the findings of Sections~\ref{sec3} and \ref{sec4} in order to identify the springs that reach plastic deformation and to estimate the time to plastic deformation in terms of the mechanical properties of networks of elastoplastic springs (Propositions~\ref{step3prop}, \ref{step4prop}, and \ref{prop5}). Section~\ref{guide_} shows the efficiency of the methodology of Section~\ref{guide} works for a particular instructive example (taken from Rachinskiy \cite{Rachinskiy}). By using the guide of Section~\ref{guide}, for the sample elastoplastic system under consideration, Section~\ref{guide_} discovers groups of indexes of the springs that are capable to reach plastic deformation and provide a sufficient condition for each of the groups to take place (Propositions~\ref{prop7}, \ref{prop7-8}, \ref{prop9}). All algebraic computations are implemented in Wolfram Mathematica which notebook is uploaded as supplementary material. We conclude Section~\ref{guide_} by remarks on the dynamics of our particular network of elastoplastic springs that Theorem~\ref{mainthm} is not capable to catch (subsection~\ref{Remarks-Dima}). Conclusions are discussed in final Section~\ref{concl}.

\vskip0.2cm

\noindent The paper includes two Appendixes. Appendix~\ref{skipped} contains proofs of not straightforward implications that we skipped proving in the main text. Appendex~\ref{appB} collects more substantial auxiliary results along with their proofs.


\section{A sufficient condition for finite-time stability}\label{sec2}

\noindent We remind the reader that the normal cone $N_C^A(y)$ to the set $C$ at a point $y\in C$ in a scalar product space $\mathcal{V}$ with the scalar product 
\begin{equation}\label{scalarproduct}
(x,y)_A=\left<x,Ay\right>,\quad{\rm where\ }A{\rm \ is\ a \ diagonal\ positive\ }m\times m{\rm -matrix,}
\end{equation} is defined as (see Bauschke-Combettes \cite[\S6.4]{Bauschke-Combettes})
$$
  N_{C}^A(y)=\left\{\begin{array}{ll}\left\{x \in \mathcal{V}:\langle x,A(\xi-y)\rangle\leqslant 0,\ {\rm for\ any }\ \xi\in {C}\right\},&\ {\rm if}\ y\in {{C}},\\
   \emptyset,& \ {\rm if}\ y\not\in {C}.
\end{array}\right.
$$

\noindent In what follows (see Bauschke-Combettes \cite[\S3.2]{Bauschke-Combettes})
\begin{eqnarray}
\|x\|^A&=&\sqrt{\left<x,Ax\right>}, \nonumber \\
{\rm proj}^A(v,F)&=&\underset{v'\in F}{\rm argmin}\|v-v'\|^A.\nonumber
\end{eqnarray}

\begin{theorem}\label{mainthm}
Let $\mathcal{V}$ be a $d$-dimensional linear subspace of $\mathbb{R}^m$ with scalar product (\ref{scalarproduct}), and  
 $c:[0,\infty)\to \mathcal{V}$ be Lipschitz continuous, and $F\subset C\subset \mathcal{V}$ be closed convex sets.
 Assume that there exists an $\eps>0$ such that condition (\ref{A0}) holds on an interval $[0,\tau_d]$ with
\begin{equation}\label{taud}
\tau_d\ge \dfrac{1}{\eps}\cdot\max_{\def\arraystretch{0.5}\begin{array}{c}_{v_1\in C}\\ _{v_2\in F}\end{array}}\|v_1-v_2\|^A.
\end{equation}
Then, every solution $y$ of (\ref{sp}) with 
$$
   C(t)=C+c(t)
$$
and any initial condition $y(0)\in C(0)$, satisfies     
 $y(\tau_d)\in F+c(\tau_d).$
\end{theorem}

\noindent What we will effectively prove is that the function
\begin{equation}\label{Vv}
   V(v)=\left<v-{\rm proj}^A(v,F),A\left(v-{\rm proj}^A(v,F)\right)\right>
\end{equation}
is a Lyapunov function for the sweeping process 
\begin{equation}\label{sp1}
  -\dot x(t)-\dot c(t)\in N_{C}^A(x(t)),
\end{equation}
which is related to (\ref{sp}) through the change of the variables (\ref{x(t)}).
Since (see Proposition~\ref{proj}) $${\rm proj}^A(v,F)+c={\rm proj}^A(v+c,F+c),$$ we have
$$
V(x(t))=\left(\|x(t)-{\rm proj}^A(x(t),F)\|^A\right)^2=\left(\|y(t)-{\rm proj}^A(y(t),F+c(t))\|^A\right)^2,
$$
for the function $x(t)$ given by (\ref{x(t)}).
Therefore, as expected, $V(x(t_1))=0$ will imply $y(t_1)\in F+c(t_1).$

\vskip0.2cm 

\noindent In what follows, $D_\xi f(u)$ is the bilateral directional derivative (Giorgi et al \cite[\S2.6]{Giorgi}, Correa-Thibault \cite {Correa-Thibault}) of $f:\mathcal{V}\to \mathcal{V}_1$ at the point $u\in \mathcal{V}$ and in the direction $\xi\in \mathcal{V},$ i.e. 
$$
D_\xi f(u)=\lim_{\tau\to 0}\dfrac{f(v+\tau\xi)-f(v)}{\tau}.
$$
Here $\mathcal{V}_1$ are finite-dimensional scalar product spaces.

\vskip0.2cm

\noindent If  the bilateral directional derivative
$D_\xi{\rm proj}^A(\cdot,F)(v)$ 
of $v\mapsto{\rm proj}^A(v,F)$ at the point $v\in \mathcal{V}$  in the direction $\xi\in \mathcal{V}$ exists, then 
the existence of $D_\xi V(v)$ and the formula 
\begin{equation}\label{V0}
   D_\xi V(v)=2\left<\xi-D_\xi{\rm proj}^A(\cdot,F)(v),A(v-{\rm proj}^A(v,F))\right>
\end{equation}
follow by observing that $$\left(\|v-{\rm proj}^A(v,F)\|^A\right)^2=\left<v-{\rm proj}^A(v,F),A(v-{\rm proj}^A(v,F))\right>,$$ see Lemma~\ref{techlemma1}. What we significantly use in the proof of Theorem~\ref{mainthm} is that any directional derivative of $v\mapsto{\rm proj}^A(v,F)$ is orthogonal to $v-{\rm proj}^A(v,F)$, so that formula (\ref{V0}) reduces to  (see Lemma~\ref{centrallemma} in Appendix~\ref{appB})
\begin{equation}\label{V1}
   D_\xi V(v)=2\left<\xi,A(v-{\rm proj}^A(v,F))\right>,
   \end{equation}
   meaning that $D_\xi V(v)$ is actually linear in $\xi$.

\begin{remark}\label{remarkblue1}
	The sweeping process possesses the property of {\em rate-independence} \cite{visintin}. Namely, if $t(\tau)$ is an increasing absolutely continuous change of time and $y(t)$ is a solution of the differential inclusion \eqref{sp} with the input $C(t)$, then $\hat y(\tau)=y(t(\tau))$ is a solution
	of the sweeping process
	\[
	-\frac{d \hat y}{d\tau}\in N^A_{\hat C(\tau)}(\hat y(\tau)), \quad \text{where} \quad \hat C(\tau)=C+\hat c(\tau),\ \  \hat c(\tau)=c(t(\tau)).
	\]
In particular, if $t(\tau)$ is taken to be the inverse of 
$$
   \tau(t)=\int_0^t |\dot c(\theta)|d\theta,
$$
then 
$$
   \hat c\hskip0.05cm'(\tau)=\dfrac{\dot c(t(\tau))}{|\dot c(t(\tau))|},
$$
so that conditions (\ref{A0}) and (\ref{taud}) of Theorem~\ref{mainthm} can be replaced by 
$$
\dot c(t)\ne 0 \ \ {\rm and} \ \	-\frac{\dot c(t)}{|\dot c(t)|}+B_\eps^A(0)\cap N_F^A(y)\subset N_{C}^A(y),\quad y\in F,\  \ {\rm a.a.}\  t\in[0,\tau_d],
$$
 and 
$$
	\int_0^{\tau_d} |\dot c(t)|\,dt \ge \dfrac{1}{\eps}\cdot\max_{\def\arraystretch{0.5}\begin{array}{c}_{v_1\in C}\\ _{v_2\in F}\end{array}}\|v_1-v_2\|^A
$$
respectively.
\end{remark}

\vskip0.2cm

\noindent {\bf Proof of Theorem~\ref{mainthm}.} Let $y(t)$ be an arbitrary solution of (\ref{sp}). For the function $x(t)$ given by (\ref{x(t)}) consider
$$
   v(t)=V(x(t)).
$$
Note, that $x(t)$ is differentiable almost everywhere on $[0,\infty)$ because $c(t)$ is Lipschitz continuous.
Since $v\mapsto{\rm proj}^A(v,F)$ is Lipschitz continuous (see e.g. Bauschke-Combettes \cite[Proposition~4.16]{Bauschke-Combettes}), the function $t\mapsto {\rm proj}^A(x(t),F)$ is differentiable almost everywhere on $[0,\infty).$ Let us fix some $t\ge 0$ such that both ${\rm proj}^A(x(t),F)$ and $x(t)$ are differentiable at $t$. Then $D_{\dot x(t)}{\rm proj}^A(\cdot,F)(x(t))$ exists (see Lemma~\ref{techlemma2} below) and by Lemma~\ref{centrallemma} we conclude
$$
   D_{\dot x(t)} V(x(t))=2\left<\dot x(t),A(x(t)-{\rm proj}^A(x(t),F))\right>.
$$
Without loss of generality we can assume that $t\ge 0$ is chosen also so that $V(x(t))$ is differentiable at $t.$ Then (see Lemma~\ref{techlemma2}),
\begin{equation}\label{dotv}
   \dot v(t)=D_{\dot x(t)} V(x(t))=2\left<\dot x(t),A\left(x(t)-{\rm proj}^A(x(t),F)\right)\right>.
\end{equation}
By the definition of  normal cone, (\ref{sp1}) implies
$$
   \left<-\dot x(t)-\dot c(t),A(\xi-x(t))\right>\le 0,\quad {\rm for \ any}\ \xi\in C.
$$
Therefore, taking $\xi={\rm proj}^A(x(t),F)$ we conclude from (\ref{dotv}) that
\begin{equation}\label{dotv1}
    \dot v(t)\le  2\left<-\dot c(t),A\left(x(t)-{\rm proj}^A(x(t),F)\right)\right>.
\end{equation}
Now we use assumption (\ref{A0}), which is equivalent to
$$
  -\dot c(t)+\eps\dfrac{\zeta}{\|\zeta\|^A}\in N_C^A(v),\quad {\rm for \ any}\ \zeta\in N_F^A(v),\ v\in F,
$$
or, using the definition of the normal cone, 
$$
  \left<-\dot c(t)+\eps\dfrac{\zeta}{\|\zeta\|^A},A(\xi-v)\right>\le 0,\quad {\rm for \ any}\ \zeta\in N_F^A(v),\ v\in F,\ \xi\in C.
$$
Therefore, letting $\xi=x(t),$ $v={\rm proj}^A(x(t),F)$, and $\zeta=x(t)-{\rm proj}^A(x(t),F),$ we get
$$
  \left<-\dot c(t)+\eps\dfrac{x(t)-{\rm proj}^A(x(t),F)}{\|x(t)-{\rm proj}^A(x(t),F)\|^A},A\left(x(t)-{\rm proj}^A(x(t),F)\right)\right>\le 0,
$$
which allows to further rewrite inequality ({\ref{dotv1}) as
$$
   \dot v(t)\le -2\eps\left<\dfrac{x(t)-{\rm proj}^A(x(t),F)}{\|x(t)-{\rm proj}^A(x(t),F)\|^A},A(x(t)-{\rm proj}^A(x(t),F))\right>=-2\eps \sqrt{v(t)}.
$$
Therefore, the Lyapunov function (\ref{Vv}) satisfies estimate (\ref{lyap2}). The proof is complete.\qed

\section{Finite-time convergence to a vertex}\label{sec3}

\noindent In this section we consider the case of ${\rm ri}(F)=\emptyset$
  or, equivalently, $M=0.$  When $M=0$, formula (\ref{F(t)}) reduces to
\begin{equation}
F=\bigcap\limits_{(\alpha,j)\in I_0} \overline{\overline{L}}(\alpha,j).\label{Fred1}
\end{equation}
In this case, of all the conditions (\ref{prop1}), (\ref{noother}), (\ref{notcoincide}), (\ref{notreduced}), (\ref{independent}), and (\ref{Ffeasible}), only conditions  
(\ref{prop1}) 
and (\ref{Ffeasible}) are needed. These two conditions take the following form:
\begin{eqnarray}
\hskip-1cm && \mbox{Condition }(\ref{prop1}):  
 \ \   \{y:y\in\overline{\overline{L}}(\alpha,j),\ (\alpha,j)\in I_0\}\mbox{ is a singleton }\{y_{*,0}\}\not=\emptyset, \label{prop1_}\\
\hskip-1cm && \mbox{Condition }(\ref{Ffeasible}):  
 \ \  y_{*,0}\in C\label{Ffeasible_}.
\end{eqnarray}

\vskip0.2cm

\begin{corollary}\label{cor1} 
Let $\mathcal{V}$ be a $d$-dimensional linear subspace of $\mathbb{R}^m$ with scalar product (\ref{scalarproduct}), and  $c:[0,\infty)\to\mathcal{V}$ be Lipschitz continuous. Assume that $C$ is given by (\ref{C(t)}) and $F$ is given by (\ref{Fred1}) with conditions (\ref{prop1_}) and (\ref{Ffeasible_}) satisfied. 
 Assume that there exists an $\eps>0$ such that  
 \begin{equation}\label{A0reduced}
  -\dot c(t)+B_\eps^A(0)\subset N_C^A(y_{*,0}),\quad t\in[0,\tau_d], \quad  \tau_d=\dfrac{1}{\eps}\cdot\max_{\def\arraystretch{0.5}\begin{array}{c}_{v\in C}\end{array}}\|v-y_{*,0}\|^A.
\end{equation}
 Then, every solution $y$ of sweeping process (\ref{sp}) with the initial condition $y(0)\in C(0)$ satisfies $y(\tau_d)=y_{*,0}+c(\tau_d).$ Furthermore, let $y_*(t)$, $t\ge \tau_d,$ be the solution of (\ref{sp}) with the initial condition $y_*(\tau_d)=y_{*,0}+c(\tau_d).$ If $c(t)$ is $T$-periodic with $T\ge\tau_d$, then $y_*$ is a globally one-period stable $T$-periodic solution of  (\ref{sp}).
\end{corollary}

\noindent We remind the reader that solution of an initial-value problem for a sweeping processes with Lipschitz continuous moving set exists, unique and features continuous dependence on initial conditions (see e.g. Kunze and Monteiro Marques \cite[Theorems 1-3]{kunze}).

\vskip0.2cm

\noindent {\bf Proof of Corollary~\ref{cor1}.} Since 
$$
   N_{\{y_{*,0}\}}^A(y_{*,0})=\mathcal{V},
$$
the inclusion of (\ref{A0}) takes the form of that of (\ref{A0reduced}). Therefore, 
by Theorem~\ref{mainthm}, $y(\tau_d)=y_{*,0}+c(\tau_d).$ Since $c(t)$ is Lipschitz continuous, sweeping process (\ref{sp}) features uniqueness of solutions and so $y(t)=y_*(t),$ $t\ge \tau_d$. Since $c(t)$ is $T$-periodic, sweeping process (\ref{sp}) admits a $T$-periodic solution (by Brouwer fixed point theorem). Therefore, $y_*$ is unique $T$-periodic solution of (\ref{sp}) (on $[\tau_d,\infty)$). Therefore, $y_*$ is the attractor of (\ref{sp}) by Massera-Krejci theorem (see Krejci \cite[Theorem
3.14]{Krejci} or Gudoshnikov-Makarenkov \cite[Theorem~4.6]{ESAIM}). The proof is complete.\qed

\vskip0.2cm

\noindent More tools about about $N_F^A(y)$ and $N_C^A(y)$ are required  to draw an applicable corollary of Theorem~\ref{mainthm} in the case where  ${\rm ri}(F)\not=\emptyset.$

\section{Finite-time convergence to a facet}\label{sec4}


\noindent 
Assume now that ${\rm ri}(F)\not=\emptyset$.
 To compute $N_C^A(y)$ and $N_F^A(y)$ for $y\in F$, we want to use the following corollary of (Rockafellar-Wets \cite[Theorem~6.46]{Rockafellar-Wets}).  Recall that ${\rm cone}\{\xi_1,...,\xi_K\}$ stays for the cone formed by vectors 
$\xi_1,...,\xi_K.$

\begin{lemma} \label{N_C} 
Let $\mathcal{V}$ be a $d$-dimensional linear subspace of $\mathbb{R}^m$ with scalar product (\ref{scalarproduct}). Consider
\begin{equation}\label{widetildeC}
   \widetilde{C}=\bigcap_{k=1}^K\left\{y\in \mathcal{V}:\left<\widetilde n_k,Ay\right>\le c_k\right\},
\end{equation}
where $\widetilde n_k\in \mathcal{V}$, $c_k\in\mathbb{R},$ $K\in\mathbb{N}.$ If $\widetilde I(y)=\left\{k\in\overline{1,K}:\left<\widetilde n_k,Ay\right>=c_k\right\}$, then 
$$
   N_{\widetilde C}^A(y)={\rm cone}\left\{\widetilde n_k:k\in \widetilde I(y)\right\}.
$$
\end{lemma}

\noindent Both, the statement of \cite[Theorem~6.46]{Rockafellar-Wets} and a proof of Lemma~\ref{N_C} are given appendix~\ref{appB}. 

\vskip0.2cm

\noindent In what follows, we will call $\widetilde n_k,\ k\in \widetilde I(y),$  the active normal vectors of the set $\widetilde C$ at a point $y$. 

\vskip0.2cm

\noindent To apply Lemma~\ref{N_C}, we rewrite $C$ and $F$ in terms of the elements of the space $\mathcal{V}.$ To this end, we consult Gudoshnikov-Makarenkov \cite[formula (27)]{PhysicaD}, which clarifies that
\begin{equation}\label{nj}\def\arraystretch{1.5}
\begin{array}{l}
   \left<e_j,A y\right>=\left<n_j,Ay\right>, \ \ n_j=\mathcal{V}_{basis}\bar e_j, \\ \bar e_j=W^{-1}\left(\begin{array}{c} R^T\\ (D^\perp)^T\end{array}\right)e_j,\ \ j\in\overline{1,m},\ y\in \mathcal{V}. 
   \end{array}
\end{equation}




\noindent To match this with the format of formula (\ref{widetildeC}), we rewrite $L(\alpha,j)$ and $\overline{\overline{L}}(\alpha,j)$ as 
\begin{eqnarray}
L(\alpha,j)&=&\left\{y\in \mathcal{V}:\left<\alpha n_j,Ay\right>\le \alpha c_j^\alpha\right\}, \label{aaa1}\\
\overline{\overline{L}}(\alpha,j)&=&\left\{y\in \mathcal{V}:\left<- n_j,Ay\right>\le- c_j^\alpha\right\}\cap
\left\{y\in \mathcal{V}:\left<n_j,Ay\right>\le  c_j^\alpha\right\}.\label{aaa2}
\end{eqnarray}
Using the representations~(\ref{aaa1})-(\ref{aaa2}), we can formulate the active normals of $C$ and $F$ at $y$ as follows:
\begin{eqnarray}
\hskip-1cm \mbox{active normals of }C\mbox{ at }y\in F&:& \{\alpha n_j:y\in\overline{\overline{L}}(\alpha,j)\}\cup\{\alpha n_j:(\alpha,j)\in I_0\},\label{line1}\\
\hskip-1cm\mbox{active normals of }F\mbox{ at }y\in F&:& \{\alpha n_j:y\in\overline{\overline{L}}(\alpha,j)\}\cup\{-n_j,n_j:(\alpha,j)\in I_0\}. \label{line2}
\end{eqnarray} 
Formula (\ref{line1}) uses condition (\ref{Ffeasible}) to make sure that the term $\{\alpha n_j:(\alpha,j)\in I_0\}$ is a part of $\{\alpha n_j:y\in\overline{\overline{L}}(\alpha,j)\}.$ 
We keep this ``redundant'' term to ease the comparison between the formulas (\ref{line1}) and (\ref{line2}). Formula (\ref{line2}) uses assumption (\ref{noother}) to claim that all vectors of $\{\alpha n_j:y\in\overline{\overline{L}}(\alpha,j)\}$ are normal vectors of $F$.

\vskip0.2cm




\noindent Using assumption (\ref{notcoincide}) we can conclude from (\ref{line1})-(\ref{line2}) that the sets of active normals are given by
\begin{eqnarray*}
 \mbox{active normals of }C\mbox{ at }y_{*,i}&:& \{\alpha n_j:(\alpha,j)\in I_i\}\cup\{\alpha n_j:(\alpha,j)\in I_0\},\ \ i\in\overline{1,M},\\
\mbox{active normals of }F\mbox{ at }y_{*,i}&:& \{\alpha n_j:(\alpha,j)\in I_i\}\cup\{-n_j,n_j:(\alpha,j)\in I_0\},\ \ i\in\overline{1,M}. 
\end{eqnarray*} 

\noindent Therefore, by Lemma~\ref{N_C}, 
\begin{equation}\label{goodformula2} \def\arraystretch{1.5}
\begin{array}{ll}
   N_C^A(y_{*,i})={\rm cone}\left\{\alpha n_j: (\alpha,j)\in I_0,\ \alpha n_j:(\alpha,j)\in I_i\right\}, & \quad i\in\overline{1,M},\\
   N_F^A(y_{*,i})={\rm cone}\left\{-n_j,n_j: (\alpha,j)\in I_0,\ \alpha n_j:(\alpha,j)\in I_i\right\}, & \quad i\in\overline{1,M}.
\end{array}
\end{equation}



\noindent Using assumption (\ref{noother}) we can specify the lists of active normals at $y\in F$  as follows:
\begin{eqnarray*}
\mbox{active normals of }C\mbox{ at }y\in F&:& \{\alpha n_j:(\alpha,j) \in\underset{i\in J_y}{\cap} I_i\}\cup\{\alpha n_j:(\alpha,j)\in I_0\},\\
\mbox{active normals of }F\mbox{ at }y\in F&:& \{\alpha n_j:(\alpha,j) \in\underset{i\in J_y}{\cap} I_i\}\cup\{-n_j,n_j:(\alpha,j)\in I_0\}.
\end{eqnarray*} 
Therefore, by (\ref{goodformula2}),
\begin{equation}\label{goodformula3}
N_C^A(y)=\bigcap_{i\in J_y} N_C^A(y_{*,i}),\quad  N_F^A(y)=\bigcap_{i\in J_y} N_F^A(y_{*,i}),\quad y\in F,
\end{equation}
where $J_y$ are given by assumption (\ref{noother}).
  
\vskip0.2cm

\noindent Provided that condition~(\ref{prop1}) holds, the boundary of the $d$-dimensional cone $N_C^A(y_{*,i})$ is the following union of $(d-1)$-dimensional cones 
\begin{equation}\label{d-1}
   \partial N_C^A(y_{*,i})=\bigcup_{(\alpha_*,j_*)\in I_0\cup I_i}{\rm cone}\left\{\alpha n_{j}:(\alpha,j)\in (I_0\cup I_i)\backslash\{(\alpha_*,j_*)\}\right\},
\end{equation}
see appendix~\ref{skipped} for a proof.
  

\vskip0.2cm

\noindent Now we use assumption (\ref{notreduced}) for the first time. This assumption allows us to conclude from (\ref{line1})-(\ref{line2}) that
\begin{eqnarray*}
\hskip-1cm \mbox{active normals of }C\mbox{ at }y\in {\rm ri}(F)&:& \ \{\alpha n_j:(\alpha,j)\in I_0\},\\
\hskip-1cm\mbox{active normals of }F\mbox{ at }y\in {\rm ri}(F)&:&\  \{-n_j,n_j:(\alpha,j)\in I_0\}. \end{eqnarray*}

\noindent Therefore, by Lemma~\ref{N_C}, 
\begin{equation}\label{goodformula1} \def\arraystretch{1.5}
\begin{array}{ll}
   N_C^A(y)={\rm cone}\left\{\alpha n_j: (\alpha,j)\in I_0\right\}, & \quad y\in{\rm ri}(F),\\
   N_F^A(y)={\rm cone}\left\{- n_j, n_j: (\alpha,j)\in I_0\right\},& \quad y\in{\rm ri}(F).
\end{array}
\end{equation}

\begin{corollary}\label{thm3} Let $\mathcal{V}$ be a $d$-dimensional linear subspace of $\mathbb{R}^m$ with scalar product (\ref{scalarproduct}), and $c:[0,\infty)\to\mathcal{V}$ be Lipschitz continuous.  Assume that $F(t)$ is given by (\ref{F(t)}) with $F$ satisfying properties (\ref{prop1})-(\ref{notcoincide}), 
and (\ref{Ffeasible}). Let $y_{*,i}$ be the vertices of $F$ given by (\ref{prop1}).
If there exists $\eps>0$ such that
\begin{eqnarray}
&&  -c'(t)+B_\eps^A(0)\cap N_F^A(y_{*,i})\subset N_C^A(y_{*,i}),\  t\in[0,\tau_d],\ i\in\overline{1,M}, \label{conc11} \\ 
 && \mbox{where}\ \tau_d=\dfrac{1}{\eps}\max\limits_{\def\arraystretch{0.5}\begin{array}{c}_{v\in C,\ i\in\overline{1,M}}\end{array}}\|v-y_{*,i}\|^A,\nonumber
\end{eqnarray}
then every solution $y$ of (\ref{sp}) with the initial condition $y(0)\in C(0)$ satisfies $y(\tau_d)\in F(\tau_d).$
\end{corollary}

\noindent {\bf Proof.} We need to prove that (\ref{conc11}) implies (\ref{A0}). Formulas (\ref{goodformula3}) and (\ref{conc11}) allow to conclude that, for any $y\in F,$
$$
   -c'(t)+B_\eps^A(0)\cap N_F^A(y)\subset -c'(t)+B_\eps^A(0)\cap N_F^A(y_{*,i})\subset N_C^A(y_{*,i}),\quad i\in J_y.
$$
Therefore,
$$
   -c'(t)+B_\eps^A(0)\cap N_F^A(y)\subset \bigcap_{i\in J_y} N_C^A(y_{*,i}),\quad y\in F,
$$
which implies (\ref{A0}) thanks to formula (\ref{goodformula3}).\qed

\vskip0.4cm

\noindent In what follows, $[W]_k$ stays for the matrix formed by first $k$ lines of the matrix $W.$

\begin{lemma}\label{lemmaforstep3} 
Let $\mathcal{V}$ be a $d$-dimensional linear subspace of $\mathbb{R}^m$ with scalar product (\ref{scalarproduct}).  Let the polyhedron $C$ and its facet $F$ be given by formulas (\ref{C(t)}) and (\ref{F(t)}). Assume that $|I_0|<d$ and
\begin{equation}\label{needlabel}
-c_1+B_\eps^A(0)\cap N_F^A(y)\subset N_C^A(y),\quad {for\ any\ }y\in{\rm ri}(F),
\end{equation}
and for some $c_1\in \mathcal{V}.$ If conditions~(\ref{prop1})-(\ref{Ffeasible}) hold, then, for any $i\in\overline{1,M}$, there exists an $\eps_i>0$ such that
$$
-c_1+B_{\eps_i}^A(0)\cap N_F^A(y_{*,i})\subset N_C^A(y_{*,i}),
$$
where 
\begin{equation}\label{formulaepsi}
\eps_i=\eps/\|\mathcal{L}_i\|^A
\end{equation}
and $\mathcal{L}_i$ is the $m\times m$-matrix given by
\begin{equation}\label{mathcalL}
\begin{array}{l}
   \mathcal{L}_i=(\{n_j,\ (\alpha,j)\in I_0\})\circ\\
   \circ\left[\left(\left(\begin{array}{c}  R^T \\ (D^\perp)^T\end{array}\hskip-0.05cm\right)(\{e_j,\ (\alpha,j)\in I_0\},\{e_j,\ (\alpha,j)\in I_i\})\right)^{-1}\right]_{|I_0|} \hskip-0.1cm\left(\begin{array}{c}  R^T \\ (D^\perp)^T\end{array}\hskip-0.05cm\right).
\end{array}
\end{equation}
\end{lemma}

\noindent {Here and in the sequel we use $''\circ''$ to denote matrix multiplication when it is broken by a line break.}

\vskip0.2cm

\noindent {\bf Proof.} Conditions (\ref{prop1}) and (\ref{independent}) imply that the vectors $\{n_j:(\alpha,j)\in I_0\cup I_i\}$ form a basis of $V$. Therefore, 
taking into account (\ref{notreduced}),
any  $\xi\in N_F^A(y_{*,i})$ is uniquely decomposable as
\begin{equation}\label{repres}
   \xi=\xi_1+\xi_2,\quad \xi_1\in N_F^A(y),\ \ \xi_2\in{\rm span}\{n_j:(\alpha,j)\in I_i\},
\end{equation}
see formulas
(\ref{goodformula2}) and (\ref{goodformula1}). Therefore, $\xi_1=\mathcal{L}_i\xi$, where $\mathcal{L}_i$ is some linear transformation (see the computation of $\mathcal{L}_i$ later in the proof) and we can estimate the norm of $\xi_1$ as 
$$
    \|\xi_1\|^A\le \|\mathcal{L}_i\|^A\cdot\|\xi\|^A.
$$
Fix $i\in\overline{1,M.}$ 
Fix an arbitrary $\xi\in N_F^A(y_{*,i})$ with $\|\xi\|^A<\eps_i$. The representation (\ref{repres}) defines $\xi_1\in N_F^A(y)$, $y\in{\rm ri}(F)$, satisfying $\|\xi_1\|^A<\eps.$ By (\ref{needlabel}),
$$
   -c_1+\xi_1\in N_C^A(y),\quad y\in{\rm ri}(F).
$$
Therefore, by formulas (\ref{goodformula2}), (\ref{goodformula1}), and now using assumptions (\ref{noother}), (\ref{notcoincide}), and (\ref{Ffeasible}),
$$
   -c_1+\xi\in N_C^A(y)+\xi_2\subset N_C^A(y_{*,i}),\quad y\in{\rm ri}(F).
$$

\noindent {\boldmath{\bf Computation of $\mathcal{L}_i$.}} Since $\{n_j:(\alpha,j)\in I_0\cup I_i\}$ is a basis of $\mathcal{V}$, we can decompose $\xi\in \mathcal{V}$ as
$$
   \xi=(\{n_j,\ (\alpha,j)\in I_0\}, \{n_j,\ (\alpha,j)\in I_i\})\left(\begin{array}{c}\zeta_1 \\ \zeta_2\end{array}\right)
$$
for some $\zeta_1\in\mathbb{R}^{|I_0|}$ and $\zeta_2\in\mathbb{R}^{|I_i|}$. On the other hand, $\xi=\mathcal{V}_{basis}v$ for some $v\in\mathbb{R}^d.$ Combining this formula and formula (\ref{nj}) for normals $n_j$, we get
$$
   \mathcal{V}_{basis}v=\mathcal{V}_{basis} W^{-1} \left(\begin{array}{c}  R^T \\ (D^\perp)^T\end{array}\right)(\{e_j,\ (\alpha,j)\in I_0\},\{e_j,\ (\alpha,j)\in I_i\})\left(\begin{array}{c}\zeta_1 \\ \zeta_2\end{array}\right)
$$
or, equivalently,
$$
   \left[\left(\begin{array}{c}  R^T \\ (D^\perp)^T\end{array}\right)(\{e_j,\ (\alpha,j)\in I_0\},\{e_j,\ (\alpha,j)\in I_i\})\right]^{-1} Wv= \left(\begin{array}{c}\zeta_1 \\ \zeta_2\end{array}\right).
$$
Therefore,
$$
   \left[\left(\left(\begin{array}{c}  R^T \\ (D^\perp)^T\end{array}\right)(\{e_j,\ (\alpha,j)\in I_0\},\{e_j,\ (\alpha,j)\in I_i\})\right)^{-1}\right]_{|I_0|}
  \circ \left(\begin{array}{c}  R^T \\ (D^\perp)^T\end{array}\right)\xi=\zeta_1,
$$
which implies (\ref{mathcalL}).
The proof of the lemma is complete.\qed

\vskip0.2cm

\begin{corollary}\label{corforstep3} Let $\mathcal{V}$ be a $d$-dimensional linear subspace of $\mathbb{R}^m$ with scalar product (\ref{scalarproduct}).  Let the polyhedron $C(t)$ and its facet $F(t)$ be given by formulas (\ref{C(t)}) and (\ref{F(t)}). Assume that ${\rm ri}(F)\not=\emptyset$ and $c'(t)=c_1$ for all $t\ge 0$. Assume that $F$ satisfies conditions (\ref{prop1})-(\ref{Ffeasible}).
If 
$$
   -c_1\in {\rm ri}(N_C^A(y)),\quad y\in {\rm ri}(F),
$$
then there exists an $\eps>0$ such that (\ref{A0}) holds on any $[0,\tau_d]$ and, in particular, the solution $y$ of sweeping process (\ref{sp}) with any initial condition $y(0)\in C(0)$ satisfies $y(\tau_d)\in F(\tau_d)$ for all sufficiently large $\tau_d>0$.
\end{corollary}

\noindent The conclusion of Corollary~\ref{corforstep3} follows by combining Lemma~\ref{lemmaforstep3} and Corollary~\ref{thm3}. The assumption on $c_1$ of Corollary~\ref{corforstep3} implies that the respective assumption of Lemma~\ref{lemmaforstep3} holds for some $\eps>0$.



\vskip0.2cm

\noindent The statement of the following remark is a part of the proof of \cite[Proposition~3.14]{ESAIM}.
\begin{remark}\label{ESAIMremark} Both $\max_{\def\arraystretch{0.5}\begin{array}{c}_{v\in C,\ i\in\overline{1,M}}\end{array}}\|v-y_{*,i}\|^A$ from Corollary~\ref{thm3} and $\max_{\def\arraystretch{0.5}\begin{array}{c}_{v\in C}\end{array}}\|v-y_{*,0}\|^A$ from Corollary~\ref{cor1} can be estimated using the following inequality
\begin{equation}\label{ESAIMestimate}
  \max_{u,v\in C}\|u-v\|^A\le \|A^{-1}c^+-A^{-1}c^-\|^A.
\end{equation}
\end{remark}
For completeness, we included a proof of formula~(\ref{ESAIMestimate}) in Appendix~\ref{skipped}.

\section{Finite-time stability of elastoplastic systems with uniaxial displacement-controlled loading}\label{mechanics}

\noindent We remind the reader that according to Moreau \cite{Moreau} a network of $m$ elastoplastic springs on $n$ nodes with 1 displacement-controlled loading is fully defined by an $m\times n$ kinematic matrix $D$ of the topology of the network, $m\times m$ matrix of stiffnesses (Hooke's coefficients) $A={\rm diag}(a_1,...,a_m)$,  an $m$-dimensional hyperrectangle $C=\prod_{j=1}^m [c_j^-,c_j^+]$ of the achievable stresses of springs (beyond which plastic deformation begins), a vector $R\in\mathbb{R}^m$ of the location of the displacement-controlled loading, and a scalar function $l(t)$ that defines the magnitude of the displacement-controlled loading. When all springs are connected (form a connected graph), we have (see Bapat \cite[Lemma 2.2]{Bapat}) 
\begin{equation}\label{rankD}
   {\rm rank}\hskip0.05cm D=n-1.
\end{equation}
We furthermore assume that
\begin{equation}\label{extra}
   m\ge n\quad{\rm and}\quad {\rm rank}(D^T R)=1.
\end{equation}


\noindent To formulate the Moreau sweeping process corresponding to the elastoplastic system $(D,A,C,R,l(t))$, we follow the 3 steps described in Gudoshnikov-Makarenkov \cite[\S5]{PhysicaD}:

\begin{itemize}

\item[1.] Find an $n\times (n-2)-$matrix $M$ of ${\rm rank}(DM)=n-2$ that solves 
$  R^T DM=0$ and use $M$ to  introduce 
$
  \mathcal{U}_{basis}=DM.
$

\item[2.] Find a matrix $\mathcal{V}_{basis}$ of $m-n+2$ linearly independent column vectors of $\mathbb{R}^m$ that solves
$  (\mathcal{U}_{basis})^T A\mathcal{V}_{basis}=0.$

\item[3.] Find an $m\times(m-n+1)-$matrix $D^\perp$ that solves 
$
  (D^\perp)^T D={\color{black}0_{(m-n+1)\times n}}
$
and such that
\begin{equation}\label{rankDperp}
{\rm rank}(D^\perp)=m-n+1. 
\end{equation}
\end{itemize}
With the new matrices introduced, the moving constraint $C(t)$ of sweeping process (\ref{sp}) corresponding to the elastoplastic system $(D,A,C,R,l(t))$ is given by 
 $$
  C(t)=\bigcap_{j=1}^m\left\{y\in \mathcal{V}:c_j^-\le\left<e_j,Ay+A\mathcal{V}_{basis}\bar L l(t)\right>\le c_j^+\right\}, 
$$
where, for each $j\in\overline{1,m},$
\begin{equation}\label{barL}
\begin{array}{l}
     \bar L=W^{-1}\left(\begin{array}{c}
      1\\ 0_{m-n+1}\end{array}\right),\ \
W=\left(\begin{array}{c}
      R^T \\
      (D^\perp)^T\end{array}\right)\mathcal{V}_{basis}
\end{array}
\end{equation}
with $e_j$ being the basis vectors of $\mathbb{R}^m$, i.e. $e_j=(\underbrace{0,...,0}_{j-1},1,0,...,0)^T.$ In this sweeping process, the variable $y(t)$ is given by
\begin{equation}\label{y(t)form}
   y(t)=A^{-1}s(t)-\mathcal{V}_{basis}\bar L l(t),
\end{equation}
where $s(t)=(s_1(t),...,s_m(t))^T$ is the vector of stresses of the springs of elastoplastic system $(D,A,C,R,l(t))$. In other words, $y(t)+\mathcal{V}_{basis}\bar L l(t)$ is the vector of elastic elongations of the springs.
It remains to observe that $C(t)$ can be rewritten in the form (\ref{C(t)}) by letting
\begin{equation}\label{c(t)0}
   \mathcal{V}=\mathcal{V}_{basis}\mathbb{R}^d,\quad d=m-n+2,\quad
     c(t)=-\mathcal{V}_{basis}\bar L l(t).
\end{equation}

\noindent The existence of $W^{-1}$ is demonstrated in Gudoshnikov-Makarenkov \cite{PhysicaD,ESAIM} for particular examples. Since this section intends to offer a general recipe, Lemma~\ref{W} in the appendix features a proof of the invertibility of $W$ in the general case.
 
\vskip0.2cm

\noindent To understand what the conclusion $y(\tau_d)\in F(\tau_d)$ of Theorem~\ref{mainthm} says about the elastoplastic system $(D,A,C,R,l(t))$, observe that formulas (\ref{y(t)form}) and (\ref{c(t)0}) imply
$$
   y(t)=A^{-1}s(t)+c(t),
$$
which upon combining with (\ref{F(t)}) gives
$$
    A^{-1}s(\tau_d)\in F.
$$
Applying Lemma~\ref{enumeration}, we conclude that the statement $y(\tau_d)\in F(\tau_d)$ is equivalent to the following property
\begin{equation}\label{enumeration1}
    s(\tau_d)\in {\rm conv}\{Ay_{*,1},...,Ay_{*,M}\}.
\end{equation}

\section{A step-by-step guide for analytic computations}\label{guide}

\noindent {\boldmath{\bf Step 1. Fix appropriate indexes $I_0$ {\it (springs  that will reach  plastic deformation)}}} 
Spot an $I_0\subset\{-1,1\}\times\overline{1,m}$ such that 
\begin{equation}\label{step1rev}
\left(\begin{array}{c}
1\\
0_{m-n+1}\end{array}\right)l'(t)\in {\rm cone}\left(
\left(\begin{array}{c}
  R^T \\
  (D^\perp)^T \end{array}\right)\left\{\alpha e_j:(\alpha,j)\in I_0\right\}\right).
\end{equation}
\begin{definition} We say that a family of indexes $I_0$ is irreducible, if $I_0$ cannot be represented in the form
\begin{equation}\label{form}
   I_0=\widetilde{I_0}\cup\{(\alpha_*,j_*)\},
\end{equation}
where $\widetilde{I_0}$ satisfies
\begin{equation}\label{step1rev_}
\left(\begin{array}{c}
1\\
0_{m-n+1}\end{array}\right)l'(t)\in {\rm cone}\left(
\left(\begin{array}{c}
  R^T \\
  (D^\perp)^T \end{array}\right)\left\{\alpha e_j:(\alpha,j)\in \widetilde{I_0}\right\}\right).
\end{equation}
\end{definition}

\noindent Proposition~\ref{irr} below explains why our results do not apply when $I_0$ is reducible. Intuitively, a vertex cannot be finite-time stable, if finite-time stability holds for the entire facet that the vertex belongs to.

\vskip0.2cm

\noindent By Corollary \ref{corrank} (see below), $I_0$ with $|I_0|=d$ always exists. However, some $I_0$ with $|I_0|=d$ may appear to be reducible, in which case  an irreducible subset of $I_0$ needs to be considered.

\begin{remark}\label{step1rem}
Relation (\ref{step1rev}) implies (see Appendix~\ref{skipped} for a proof)
\begin{equation}\label{step1conclusion}
   -c'(t)\in{\rm cone}\left\{\alpha n_j: (\alpha,j)\in I_0\right\}.
\end{equation}
\end{remark}








\vskip0.2cm

\noindent {\boldmath{\bf Step 2. Fix appropriate indexes $I_i$}} {\it (springs that may reach plastic deformation and that affect the convergence of springs $I_0$ to plastic deformation).} {\boldmath{\bf  Skip this step, if $|I_0|=d$.}} We will consider the simplest possible way to design $I_i$ which ensures  that $F\not=\emptyset$ and which satisfies the assumptions (\ref{F(t)})-(\ref{Ffeasible}). This simplest way utilizes the minimal possible number $d-|I_0|$ of springs. The conditions to be imposed on the remaining $m-d$ springs will ensure that those $m-d$ springs don't affect the convergence of the stress vector to $F(t)$ and, in particular, don't undergo plastic deformation when close to $F(t)$.

\vskip0.2cm

\noindent Find some $I_1$ such that 
\begin{equation}\label{+}
|I_0|+|I_1|=d
\end{equation} and such that 
\begin{equation}\label{conditionstep2}
{\rm rank}\left(
\left(\begin{array}{c}
  R^T \\
  (D^\perp)^T \end{array}\right)(\left\{\alpha e_j:(\alpha,j)\in I_0\cup I_1\right\})\right)=d.
\end{equation} Based on $I_1$ we can obtain more vertexes $I_i$ by changing the elements of $I_1$ from $(\alpha,n)$ to $(-\alpha,n).$ Let $I_i,$ $i\in\overline{1,M},$ 
$$
M=2^{d-|I_0|}
$$ 
be all different families of indexes obtained through this process, i.e.
\begin{equation}\label{different}
   I_{i_1}\not=I_{i_2},\quad i_1\not=i_2,\ 
   \ i_1,i_2\in \overline{1,M}. 
\end{equation}  

\vskip0.0cm

\noindent Use $I_0$ and $I_i$, $i\in\overline{1,M},$ to define $F$ by formula (\ref{F(t)}). 

\vskip0.2cm

\noindent {\boldmath{\bf Step 3. Compute the vertexes of $F$ and impose conditions ensuring feasibility of $F.$}} Depending on whether $|I_0|=d$ or $|I_0|<d,$ compute $y_{*,0}$ or $y_{*,i},$ $i\in \overline{1,M}$, using the formula (see Appendix~\ref{skipped} for a proof)
\begin{equation}\label{y*iformula}
   y_{*,i}=\mathcal{V}_{basis}\left(
   \left(\left\{e_j,(\alpha,j)\in I_0\cup I_i\right\}\right)^T      
    A \mathcal{V}_{basis}\right)^{-1}\left(\left\{c_j^\alpha,(\alpha,j)\in I_0\cup I_i\right\}\right)^T.
\end{equation}
The feasibility condition (\ref{Ffeasible}) will hold, if
\begin{equation}\label{y*ifeasible}
\begin{array}{ll}
|I_0|<d: &\ \    c_j^-<\left<e_j,Ay_{*,i}\right>< c_j^+,\ \  i\in\overline{1,M},\ (\alpha,j)\not\in I_0\cup I_1\cup ... \cup I_M,\\
|I_0|=d: &\ \   c_j^-<\left<e_j,Ay_{*,0}\right>< c_j^+,\ \  (\alpha,j)\not\in I_0.
\end{array}
\end{equation}

\vskip0.2cm

\noindent Assumption (\ref{notcoincide}) concerning non-coincidence of the vertices will hold if
\begin{equation}\label{different1}
\begin{array}{ll}
|I_0|<d: &\ \      c_j^-<c_j^+,\ \ \mbox{for all } (\alpha,j)\in I_i,\ i\in\overline{1,M},\\
|I_0|=d: &\ \   \mbox{does not apply}.
\end{array}
\end{equation}

\noindent We will say that relation (\ref{step1rev}) holds in a strict sense, if, on top of (\ref{step1rev}), the following property is satisfied:
\begin{equation}\label{step1rev-strict}
\left(\begin{array}{c}
1\\
0_{m-n+1}\end{array}\right)l'(t)\not\in {\rm rb}\left({\rm cone}\left(
\left(\begin{array}{c}
  R^T \\
  (D^\perp)^T \end{array}\right)\left\{\alpha e_j:(\alpha,j)\in I_0\right\}\right)\right).
\end{equation}

\noindent 
 With the moving constraint $C(t)$ introduced in Section~\ref{mechanics} and with the facet $F$ introduced in Steps~1-3, the Corollaries~\ref{cor1} and \ref{corforstep3} lead to the following qualitative description of the asymptotic behavior of elastoplastic system 
$(D,A,C,R,l(t))$ and associated sweeping process (\ref{sp}).

\begin{proposition} {\bf (Conclusion of Steps 1-3).} \label{step3prop} If $l'(t)=const,$ if relation (\ref{step1rev}) holds in a strict sense, and if properties (\ref{y*ifeasible}) and (\ref{different1}) hold, then there exists an $\eps>0$ such that condition (\ref{A0}) is satisfied on any $[0,\tau_d]$ and, in particular, the solution $y$ of sweeping process (\ref{sp}) with the initial condition $y(0)\in C(0)$ satisfies $y(\tau_d)\in F(\tau_d)$ for all sufficiently large $\tau_d>0$. Accordingly, the stress vector $s(t)$ of the elastoplastic system 
$(D,A,C,R,l(t))$ satisfies (\ref{enumeration1}) 
for all sufficiently large $\tau_d$.

\end{proposition}

\noindent When $|I_0|=d$, the statement of Proposition~\ref{step3prop} follows from Corollary~\ref{cor1} almost directly. Assumption~(\ref{prop1_}) holds because $I_0$ is irreducible. Assumption~(\ref{Ffeasible_}) is satisfied by (\ref{y*ifeasible}). Conditions (\ref{step1rev}) and (\ref{step1rev-strict}) ensure the existence of $\eps>0$ for which (\ref{A0reduced}) holds for any $t\ge 0.$

\vskip0.2cm

\noindent Considering $|I_0|<d$ and deriving the statement of Proposition~\ref{step3prop} from Corollary~\ref{corforstep3} requires establishing validity of assumptions (\ref{prop1}), (\ref{noother}), (\ref{notcoincide}), (\ref{notreduced}), (\ref{independent}), and (\ref{Ffeasible}). Property (\ref{prop1}) follows from (\ref{conditionstep2}). Property (\ref{notcoincide}) follows from (\ref{different}) and (\ref{different1}).
Property (\ref{Ffeasible}) follows from (\ref{y*ifeasible}).
 Property  (\ref{independent}) coincides with (\ref{+}). Verifying conditions (\ref{noother}) and (\ref{notreduced}) is less straightforward. This is done in the two propositions that follow below.

\begin{proposition}\label{propnoother} Assume $M\ge 1.$ Let $F$ be the facet defined in Step~2. If (\ref{y*ifeasible}) holds then (\ref{noother}) holds as well. In other words, (\ref{y*ifeasible}) implies that
$$
\begin{array}{c}
  \{(\alpha,j):y\in\overline{\overline{L}}(\alpha,j)\}=I_0\cup\bigcap\limits_{i\in I} I_i,\ \ \mbox{for some }I\subset\overline{1,M}\ \mbox{or for }I=\emptyset.
  \end{array}
$$
\end{proposition}

\noindent {\bf Proof.} Let $y\in F$ and let $I_*=\{(\alpha,j):y\in\overline{\overline{L}}(\alpha,j)\}.$  By condition (\ref{y*ifeasible}),
$I_*=I_0\cup I_{**},$ where $I_{**}\subset I_1\cup\ldots \cup I_M.$
By construction,
\begin{equation}\label{expr1}
   I_i=\{(\pm,n_{j_1}),\ldots,(\pm,n_{j_{d-|I_0|}})\},
\end{equation}
where different $i$ correspond to different choices of $''+''$ and $''-''$ in each symbol $''\pm''$. 
Therefore, $I_{**}$ can either be an empty set or a set of the form
\begin{equation}\label{expr2}
   I_{**}=\{(\alpha_1^{**},n_{j_1^{**}}),\ldots,(\alpha_{d-|I_0|}^{**},n_{j_{d-|I_0|}^{**}})\},
\end{equation}
where
$\{j_1^{**},\ldots,j_{d-|I_0|}^{**}\}\subset\{j_1,\ldots,j_{d-|I_0|}\}.$ If $I_{**}=\emptyset$, then the proof is complete. So, from now on we assume that $I_{**}\not=\emptyset.$ 

\vskip0.2cm

\noindent From expressions (\ref{expr1}) and (\ref{expr2}) we see that $I_{**}\subset I_i$ for at least one index $i_0\in\overline{1,M}$. 
Define $I$ as
$$
   I=\{i\in\overline{1,M}: I_{**}\subset I_i\}.
$$
Therefore 
\begin{equation}\label{turns}
\begin{array}{c}   I_{**}\subset \bigcap\limits_{i\in I} I_i.\end{array}
\end{equation}
Since the elements of $I_i\backslash I_{**}$ are obtained from the elements of $I_{i_0}\backslash I_{**}$ by taking all possible replacements of $(\alpha,n)$ by $(-\alpha,n)$, we have 
$$
\begin{array}{c}   \bigcap\limits_{i\in I} (I_i\backslash I_{**})=\emptyset.\end{array}
$$
Therefore, $\bigcap\limits_{i\in I} I_i \subset I_{**}$, and so (\ref{turns}) turns into equality.\qed

\begin{lemma}\label{riF} Assume that the facet $F$ is given by (\ref{F(t)}). Assume that conditions (\ref{prop1}) and (\ref{independent}) hold. If
\begin{equation}
 \begin{array}{l} \mbox{there exists }\bar y\in F\mbox{ such that}\  \left<ae_j,A\bar y\right><\alpha c_j^\alpha, \ (\alpha,j)\in I_i,\ i\in\overline{1,M},\end{array} \label{ass2}
\end{equation}
then condition (\ref{notreduced}) is satisfied.
\end{lemma}

\noindent {\bf Proof.} {\it Part 1.} If
 (\ref{ass2}) holds then there exists a full-dimensional ball $B_\delta(\bar y)$ in $V$ such that
$$
   B_\delta(\bar y)\subset L(\alpha,j),\quad (\alpha,j)\in I_i,\ i\in\overline{1,M}.
$$
Therefore,
$$
\begin{array}{c}  {\rm aff}(F)\supset {\rm aff}\left(\bigcap\limits_{(\alpha,j)\in I_0} \overline{\overline{L}}(\alpha,j)\cap B_\delta(\bar y)\right)={\rm aff}\left(\bigcap\limits_{(\alpha,j)\in I_0} \overline{\overline{L}}(\alpha,j)\right),\end{array}
$$
where ${\rm aff}(A)$ is the affine hull of set $A.$ 
On the other hand, directly from the definition of $F$,
$$
\begin{array}{c}  {\rm aff}(F)\subset  {\rm aff}\left(\bigcap\limits_{(\alpha,j)\in I_0} \overline{\overline{L}}(\alpha,j)\right).\end{array}
$$
So we conclude that
\begin{equation}\label{weconclude}\begin{array}{c}
  {\rm aff}(F)=   {\rm aff}\left(\bigcap\limits_{(\alpha,j)\in I_0} \overline{\overline{L}}(\alpha,j)\right).\end{array}
\end{equation}

\noindent {\it Part 2.} Consider $y\in F$ and assume that $y\in\overline{\overline{L}}(\alpha_*,j_*)$ for some $(\alpha_*,j_*)\in I_{i_*}$ and some $i_*\in\overline{1,M}.$ 
By properties~(\ref{prop1}) and (\ref{independent}), the subspace (\ref{weconclude}) intersects the subspace $\overline{\overline{L}}(\alpha_*,j_*)$ transversally. 
Therefore, if we consider a ball of space (\ref{weconclude}) centered at $y$, then part of this ball will lie outside of $L(\alpha_*,j_*).$ Therefore, $y\in{\rm rb}(F)$. Therefore, if $y\in{\rm ri}(F)$ then $y\not\in\overline{\overline{L}}(\alpha,j),$ $(\alpha,j)\in I_i$, $i\in\overline{1,M}$, which completes the proof. \qed

\begin{proposition} Assume $M\ge 1.$ Let $F$ be the facet defined in Step~2. If (\ref{different1}) holds then (\ref{ass2}) holds as well.
\end{proposition}

\noindent {\bf Proof.} We will construct the required $\bar y$ as the solution of the following system of $d$ algebraic equations:
$$
\left\{\begin{array}{ll}
\left<e_j,A\bar y\right>=c_j^\alpha, &\ \ (\alpha,j)\in I_0,\\
\left<e_j,A\bar y\right>=\dfrac{c^-_j+c_j^+}{2},&\ \ (\alpha,j)\in I_i,\ i\in\overline{1,M}.
\end{array}\right.
$$
As in the proof of formula (\ref{y*iformula}), this system of equations admits a unique solution $\bar y$ because $F$ satisfies assumptions (\ref{prop1}) and (\ref{independent}). Condition (\ref{different1}) implies that
$$
   c_j^-<\dfrac{c_j^-+c_j^+}{2}<c_j^+.
$$
Therefore,
by construction, 
\begin{eqnarray*}
&& \bar y\in \bigcap\limits_{(\alpha,j)\in I_0} \overline{\overline{L}}(\alpha,j),\\
&& c_j^-<\left<e_j,A\bar y\right><c_j^+,\quad (\alpha,j)\in I_i,\ i\in\overline{1,M},
\end{eqnarray*}
which yields (\ref{ass2}).
\qed

\vskip0.5cm

\noindent One has to proceed to Steps 4 and 5, if an estimate for $\tau_d$ is of interest.

\vskip0.2cm

\noindent {\boldmath{\bf Step 4. Compute $\eps_0.$}} Our next argument will be based on application of Corollary~\ref{cor1} (when $|I_0|=d$) and Corollary~\ref{thm3} in combination with Lemma~\ref{lemmaforstep3}. This step is devoted to finding $\eps$ for which the respective assumptions~(\ref{A0reduced}) and (\ref{needlabel}) hold. Assumptions (\ref{A0reduced}) and (\ref{needlabel}) require computing the distance from $-c'(t)$ to the boundary of cone $N_C^A(y)$ at the point 
$F$ when $F$ is a singleton and at the points of ${\rm ri}(F)$ when ${\rm ri}(F)\not=\emptyset.$ In either case, the required boundary is 
$\partial\hskip0.05cm {\rm cone}\left\{\alpha n_j:(\alpha,j)\in I_0\right\}$.

\vskip0.2cm

\noindent Using formula (\ref{d-1}), we compute 
\begin{equation}\label{newnew} 
 \def\arraystretch{1.5}
\begin{array}{l}
\eps_0(t)={\rm dist}^A\left(-c'(t),\partial {\rm cone}\left\{\alpha n_j:(\alpha,j)\in I_0\right\}\right)=\\
 =\min\limits_{(\alpha_*,j_*)\in I_0}{\rm dist}^A \left(
-c'(t),{\rm cone} \left\{\alpha n_{j}:(\alpha,j)\in I_0\backslash\{(\alpha_*,j_*)\}\right\}\right).
\end{array}
\end{equation} 

\noindent The following lemma can be used to compute the distances from $-c'(t)$ to the required cones (see Appendix~\ref{appB} for a proof of the lemma).

\begin{lemma} \label{computedistance} Assume that $\{n_{i_1},...,n_{i_k}\}$ is a linearly independent subset of vectors $\{n_1,...,n_m\}$. Introduce
$\mathcal{N}=\left(n_{i_1} \ ... \ n_{i_k}\right)$. Then the matrix 
$\mathcal{N}^T A \mathcal{N}$ is invertible and, for any $c'\in\mathbb{R}^m$,  
\begin{eqnarray}
&&   {\rm dist}^A\left(-c', {\rm span}\left\{n_{i_1},...,n_{i_k}\right\}\right)=\left\|-c'-{\rm proj}^A\left(-c',{\rm span}\left\{n_{i_1},...,n_{i_k}\right\}\right)\right\|^A,\nonumber\\
 &&  {\rm proj}^A\left(-c',{\rm span}\left\{n_{i_1},...,n_{i_k}\right\}\right)=-\mathcal{N}\left[\mathcal{N}^T A \mathcal{N}\right]^{-1}\mathcal{N}^T A c'.\label{distcone}
\end{eqnarray}
\end{lemma}

\noindent Based on Lemma~\ref{computedistance} we can rewrite formula (\ref{newnew}) as
$$
\bar \eps_0(t)=\min\limits_{(\alpha_*,j_*)\in I_0}\left\|-c'(t)
 -{\rm proj}^A\left(-c'(t),{\rm span}\left\{\alpha n_{j}:(\alpha,j)\in I_0\backslash\{(\alpha_*,j_*)\}\right\}\right)\right\|^A,
$$

\vskip-0.9cm

\begin{equation}\label{newnewnew}
\begin{array}{l}
 {\rm proj}^A\left(-c'(t),{\rm span}\left\{\alpha n_{j}:(\alpha,j)\in I_0\backslash\{(\alpha_*,j_*)\}\right\}\right)=\\
 =-(\{n_j,\ (\alpha,j)\in I_0\backslash\{(\alpha_*,j_*)\}\})\circ\\
\circ\left[(\{n_j,\ (\alpha,j)\in I_0\backslash\{(\alpha_*,j_*)\}\})^T A (\{n_j,\ (\alpha,j)\in I_0\backslash\{(\alpha_*,j_*)\}\})\right]^{-1}\circ\\
\circ (\{n_j,\ (\alpha,j)\in I_0\backslash\{(\alpha_*,j_*)\}\})^T Ac'.
\end{array}
\end{equation}


\noindent Choose $\eps_0>0$ such that $\eps_0\le \bar\eps_0(t)$ for all $t\in[0,\tau_d].$ Corollary~\ref{cor1}, Remark~\ref{ESAIMremark} and formula (\ref{enumeration1}) then lead to the following conclusion.

\vskip0.2cm

\begin{proposition} {\bf (Conclusion of Steps 1-4).} \label{step4prop}Assume that $|I_0|=d$, i.e. $F=\{y_{*,0}\}$. Assume that conditions (\ref{step1rev}) and (\ref{y*ifeasible}) hold on $[0,\tau_d]$. If
$$
  \tau_d\ge \dfrac{1}{\eps_0} \cdot\|A^{-1}c^+-A^{-1}c^-\|^A,
$$
then condition (\ref{A0}) is satisfied on $[0,\tau_d]$ and, in particular, the solutions $y$ of sweeping process (\ref{sp}) with any initial conditions $y(0)\in C(0)$ satisfy $y(\tau_d)=y_{*,0}+c(\tau_d)$. Accordingly, the stress vector $s(t)$ of the elastoplastic system 
$(D,A,C,R,l(t))$ satisfies 
$
    s(\tau_d)= Ay_{*,0}
$
regardless of the initial value $s(0).$ 
Therefore, if $c(t)$ is $T$-periodic with $T\ge\tau_d$, then the solution $y_*(t)$ with the initial condition $y_*(\tau_d)=y_{*,0}+c(\tau_d)$ is a  one-period stable $T$-periodic solution of (\ref{sp}) and the stress-vector  $s(t)$ of the elastoplastic system $(D,A,C,R,l(t))$ exhibits a unique $T$-periodic behavior beginning the time $\tau_d.$
\end{proposition}


\vskip0.0cm

\noindent One more step is required to produce an estimate for $\tau_d$ when $|I_0|<d.$

\vskip0.2cm

\noindent {\boldmath{\bf Step 5. Compute $\sigma_i.$}} Having found $\eps_0$ for which (\ref{needlabel}) holds, 
we can now use 
 Lemma~\ref{lemmaforstep3} to compute $\eps$ for which assumption~(\ref{conc11}) of Corollary~\ref{thm3} is satisfied.   Specifically, formula (\ref{formulaepsi}) of Lemma~\ref{lemmaforstep3} implies that the required $\eps$ is given by
 $$
    \eps=\eps_0\min\limits_{i\in\overline{1,M}}(1/\|\mathcal{L}_i\|^A)=\eps_0\max\limits_{i\in\overline{1,M}}\|\mathcal{L}_i\|^A,
 $$
 where the $m\times m$ matrices $\mathcal{L}_i$ are given by (\ref{mathcalL}).
 
 \vskip0.2cm
 
 \noindent In order to compute $\|\mathcal{L}_i\|^A$, we first observe that
 \begin{eqnarray*}
    \|\mathcal{L}_i\xi\|^A&=&\sqrt{\left<\mathcal{L}_i\xi,A\mathcal{L}_i\xi\right>}=\sqrt{\left<\sqrt{A}\mathcal{L}_i\xi,\sqrt{A}\mathcal{L}_i\xi\right>}=\|\sqrt{A}\mathcal{L}_i\xi\|=\\
    &=&\|\sqrt{A}\mathcal{L}_i\sqrt{A^{-1}}\sqrt{A}\xi\|\le \|\sqrt{A}\mathcal{L}_i\sqrt{A^{-1}}\|\cdot \|\sqrt{A}\xi\|=\|\sqrt{A}\mathcal{L}_i\sqrt{A^{-1}}\|\cdot \|\xi\|^A.
 \end{eqnarray*}
 Therefore,
 $$
    \|\mathcal{L}_i\|^A\le \|\sqrt{A}\mathcal{L}_i\sqrt{A^{-1}}\|.
 $$
 But, based on e.g. Friedberg et al \cite[\S6.10, Corollary~1]{Friedberg}, 
$$\|\sqrt{A}\mathcal{L}_i\sqrt{A^{-1}}\|=\sqrt{\sigma_i},$$ where 
\begin{equation}\label{eigenvalue}
\sigma_i\mbox{ is the largest eigenvalue of  matrix }\left(\sqrt{A}\mathcal{L}_i\sqrt{A^{-1}}\right)^T\sqrt{A}\mathcal{L}_i\sqrt{A^{-1}}.
\end{equation} 

\vskip0.2cm

\noindent Therefore, $\eps_i$ can be computed as
$$
  \eps_{i}=\eps_0 \hskip0.07cm /\max\limits_{i\in\overline{1,M}}\sqrt{\sigma_i}.
$$



\noindent Corollary~\ref{thm3}, Remark~\ref{ESAIMremark}, and formula (\ref{enumeration1}) can now be summarized as follows. 

\begin{proposition} {\bf (Conclusion of Steps 1-5).} \label{prop5} Assume that $|I_0|<d$. Assume that conditions (\ref{step1rev}) and (\ref{y*ifeasible}) hold on $[0,\tau_d]$. If
$$
  \tau_d\ge \dfrac{\max\{\sqrt{\sigma_1},\ldots,\sqrt{\sigma_M}\}}{\eps_{0}} \cdot\|A^{-1}c^+-A^{-1}c^-\|^A,
$$
then condition (\ref{A0}) is satisfied on $[0,\tau_d]$ and, in particular, the solutions $y$ of sweeping process (\ref{sp}) with any initial conditions $y(0)\in C(0)$ satisfy $y(\tau_d)\in F(\tau_d)$.  Accordingly, the stress vector $s(t)$ of the elastoplastic system 
$(D,A,C,R,l(t))$ satisfies 
$
    s(\tau_d)\in {\rm conv}\{Ay_{*,1},...,Ay_{*,M}\}.
$
\end{proposition}


\noindent We remind the reader that inclusion (\ref{step1rev}) is called strict, if (\ref{step1rev-strict}) holds.

\begin{proposition}\label{irr} If $I_0$ is reducible, then inclusion (\ref{step1rev}) is never strict and, in particular, $\eps_0(t)$ given by formula (\ref{newnew}) is necessarily zero.
\end{proposition}

\noindent {\bf Proof.} By definition, $I_0$ is representable in the form (\ref{form}). Therefore, as in the proof of formula (\ref{d-1}), we can conclude that
$$
  {\rm cone}\{\alpha n_j:(\alpha,j)\in \widetilde{I_0}\}\subset {\rm rb}\left({\rm cone}\left\{\alpha n_j:(\alpha,j)\in I_0\right\}\right).
$$
Hence, by (\ref{step1rev_}), 
$$
-c'(t)\in{\rm rb}\left( {\rm cone}\hskip0.05cm\left\{\alpha n_j:(\alpha,j)\in I_0\right\}\right).
$$
Therefore, inclusion (\ref{step1rev}) is not strict (we use Remark~\ref{step1rem} again) and $\eps_0(t)$ given by (\ref{newnew}) vanishes.\qed

\section{Application to a system of elastoplastic springs}\label{guide_}

\begin{table}[] 
\centering
\begin{tabular}{|l|l|l|l|}
\hline
\multicolumn{2}{|l|}{Step 1}            &                                     & Step 2                                          \\ \hline  \multirow{11}{*}{$|I_0|=2$  }  
                   &                    & \cellcolor[HTML]{C0C0C0}scenario 1 & $\begin{array}{c} I_1=\{(-,3)\}\\ I_2=\{(+,3)\}\end{array}$ \cellcolor[HTML]{C0C0C0}                        \\ \cline{3-4} 
                   &         $I_0=\{(+,1),(+,2)\}$           & scenario 2                         &  $\begin{array}{c} I_1=\{(-,4)\}\\ I_2=\{(+,4)\}\end{array} $                                               \\ \cline{3-4} 
                   & \multirow{-3}{*}{} & \cellcolor[HTML]{C0C0C0}scenario 3 & $\begin{array}{c} I_1=\{(-,5)\}\\  I_2=\{(+,5)\}\end{array} $ \cellcolor[HTML]{C0C0C0}                        \\ \cline{2-4} 
                   &                    & scenario 4                         & $\begin{array}{c} I_1=\{(-,1)\}\\  I_2=\{(+,1)\}\end{array}  $                                                 \\ \cline{3-4} 
                   &         $I_0=\{(+,4),(+,5)\}$           & \cellcolor[HTML]{C0C0C0}scenario 5 &  $\begin{array}{c} I_1=\{(-,2)\}\\  I_2=\{(+,2)\}\end{array}$ \cellcolor[HTML]{C0C0C0}                        \\ \cline{3-4} 
\multirow{-6}{*}{} & \multirow{-3}{*}{} & scenario 6                         &                                               $\begin{array}{c} I_1=\{(-,3)\}\\  I_2=\{(+,3)\}\end{array} $   \\ \hline \multirow{2}{*}{$|I_0|=3$  }   
                   &        $I_0=\{(+,1),(-,3),(+,5)\}$            & \cellcolor[HTML]{C0C0C0}scenario 7 & \cellcolor[HTML]{343434}{\color[HTML]{343434} } \\ \cline{2-4} 
\multirow{-2}{*}{} &        $I_0=\{(+,2),(+,3),(+,4)\}$            & scenario 8                         & \cellcolor[HTML]{343434}{\color[HTML]{343434} } \\ \hline
\end{tabular} 
\caption{A list of possible scenarios according to which the elastoplastic system shown in Fig.~\ref{fig1} stabilizes under an  increasing or decreasing input. The notations of the table are introduced and explained in Steps 1 and 2 of Section~\ref{guide} and computed for model of Fig.~\ref{fig1} in Steps 1 and 2 of Section~\ref{guide_}.}
\label{table1}
\end{table}

\begin{table}[]
\centering
\begin{tabular}{|l|l|l|}
\hline
                                    & \multicolumn{2}{c|}{Step 3}                                                                  \\ \cline{2-3} 
\multirow{-2}{*}{}                  &                      
$\begin{array}{c}
\mbox{\color{black} the endpoints of the line-segment}\\
  {\color{black}[
 Ay_{*,1},Ay_{*,2}]}
\\
 \mbox{\color{black}of terminal stresses of the springs}
\end{array}$                         & \multicolumn{1}{c|}{feasibility condition} \\ \hline
\rowcolor[HTML]{C0C0C0} 
\cellcolor[HTML]{C0C0C0}scenario 1 &     $\begin{array}{l} Ay_{*,1}=\left(\begin{array}{c} c_1^+ \\ c_2^+ \\ c_3^- \\ c_1^+-c_3^- \\ c_2^++c_3^-\end{array}\right) \\ Ay_{*,2}=\left(\begin{array}{c} c_1^+ \\ c_2^+ \\ c_3^+ \\ c_1^+-c_3^+ \\ c_2^++c_3^+\end{array}\right)\end{array}$                                                                                       &                                           
$ \def\arraystretch{1.5}
\begin{array}{l}
c_4^-\le c_1^+-c_3^-\le c_4^+\\
c_5^-\le c_2^++c_3^-\le c_5^+\\
\\ 
c_4^-\le c_1^+-c_3^+\le c_4^+\\
c_5^-\le c_2^++c_3^+\le c_5^+
\end{array}$
 \\ \hline
scenario 2                         &     $\begin{array}{l}Ay_{*,1}=\left(\begin{array}{c} c_1^+ \\ c_2^+ \\ c_1^+ -c_4^- \\ c_4^- \\ c_1^++c_2^+-c_4^-\end{array}\right) \\ Ay_{*,2}=\left(\begin{array}{c} c_1^+ \\ c_2^+ \\ c_1^+ -c_4^+\\ c_4^+\\ c_1^+ + c_2^+-c_4^+\end{array}\right)\end{array}$                                               &                                           
$ \def\arraystretch{1.5}
\begin{array}{rcl}
c_3^-\le &c_1^+-c_4^-&\le c_3^+\\
c_5^-\le &c_1^++c_2^+-c_4^-&\le c_5^+\\ \\
c_3^-\le &c_1^+-c_4^+&\le c_3^+\\
c_5^-\le &c_1^++c_2^+-c_4^+&\le c_5^+
\end{array}$
 \\ \hline
\rowcolor[HTML]{C0C0C0} 
\cellcolor[HTML]{C0C0C0}scenario 3 &     $\begin{array}{l}Ay_{*,1}=\left(\begin{array}{c} c_1^+ \\ c_2^+ \\ -c_2^+ + c_5^- \\ c_1^++ c_2^+-c_5^- \\ c_5^-\end{array}\right) \\ Ay_{*,2}=\left(\begin{array}{c} c_1^+ \\ c_2^+ \\ -c_2^+ +c_5^+\\ c_1^++c_2^+-c_5^+\\ c_5^+\end{array}\right)\end{array}$                                              &                                           
$ \def\arraystretch{1.5}
\begin{array}{rcl}
c_3^-\le &-c_2^++c_5^-&\le c_3^+\\
c_4^-\le &c_1^++c_2^+-c_5^-&\le c_4^+\\ \\
c_3^-\le &-c_2^++c_5^+&\le c_3^+\\
c_4^-\le &c_1^++c_2^+-c_5^+&\le c_4^+
\end{array}$
\\ \hline
scenario 4                         &         $\begin{array}{l}Ay_{*,1}=\left(\begin{array}{c} c_1^- \\ -c_1^-+c_4^++c_5^+ \\ c_1^--c_4^+ \\ c_4^+ \\ c_5^+\end{array}\right) \\ 
Ay_{*,2}=\left(\begin{array}{c} c_1^+ \\ -c_1^++c_4^++c_5^+ \\ c_1^+-c_4^+\\ c_4^+\\ c_5^+\end{array}\right)\end{array}$                                          &                                         
$ \def\arraystretch{1.5}
\begin{array}{rcl}
c_2^-\le &-c_1^-+c_4^++c_5^+&\le c_2^+\\
c_3^-\le &c_1^--c_4^+&\le c_3^+\\ \\
c_2^-\le &-c_1^++c_4^++c_5^+&\le c_2^+\\
c_3^-\le &c_1^+-c_4^+&\le c_3^+
\end{array}$
   \\ \hline
\end{tabular} 
\caption{Vertices of the attractors and the corresponding feasibility conditions for the first 4 scenarios of Table~\ref{table1}. 
The notation of the table are introduced and explained in Step~3 of Section~\ref{guide} and computed for model of Fig.~\ref{fig1} in Step~3 of Section~\ref{guide_}.}
\label{table2}
\end{table}

\begin{table}[]
\begin{tabular}{|l|l|l|}
\hline
                                    & \multicolumn{2}{c|}{Step 3}                                                                  \\ \cline{2-3} 
\multirow{-2}{*}{}                  &                  $\begin{array}{c}
\mbox{\color{black} the endpoints of the line-segment}\\
 \ {\color{black}[
 Ay_{*,1},Ay_{*,2}]}
\\
 \mbox{\color{black}of terminal stresses of the springs}\\
 \mbox{\color{black}or the terminal stress of the springs}\\
{\color{black} A y_{*,0}}
\end{array}$                               & \multicolumn{1}{c|}{feasibility condition} \\ \hline
\rowcolor[HTML]{C0C0C0} 
\cellcolor[HTML]{C0C0C0}scenario 5 &  $\begin{array}{l}Ay_{*,1}=\left(\begin{array}{c} -c_2^-+c_4^++c_5^+ \\ c_2^- \\ -c_2^-+c_5^+ \\ c_4^+ \\ c_5^+\end{array}\right) \\ 
Ay_{*,2}=\left(\begin{array}{c} -c_2^++c_4^++c_5^+ \\ c_2^+ \\ -c_2^++c_5^+ \\ c_4^+ \\ c_5^+\end{array}\right)\end{array}$                                               &                                           $ \def\arraystretch{1.5}
\begin{array}{rcl}
c_1^-\le &-c_1^-+c_4^++c_5^+&\le c_1^+\\
c_3^-\le &c_1^--c_4^+&\le c_3^+\\ \\
c_1^-\le &-c_1^++c_4^++c_5^+&\le c_1^+\\
c_3^-\le &c_1^+-c_4^+&\le c_3^+
\end{array}$
 \\ \hline
scenario 6                         &                    $\begin{array}{l}Ay_{*,1}=\left(\begin{array}{c} c_3^-+c_4^+ \\ -c_3^-+c_5^+ \\ c_3^- \\ c_4^+ \\ c_5^+\end{array}\right) \\ 
Ay_{*,2}=\left(\begin{array}{c} c_3^++c_4^+ \\ -c_3^++c^+_5 \\ c_3^+ \\ c_4^+ \\ c_5^+\end{array}\right)\end{array}$                                               &                                           $ \def\arraystretch{1.5}
\begin{array}{rcl}
c_1^-\le &c_3^-+c_4^+&\le c_1^+\\
c_2^-\le &-c_3^-+c_5^+&\le c_2^+\\ \\
c_1^-\le &c_3^++c_4^+&\le c_1^+\\
c_2^-\le &-c_3^++c_5^+&\le c_2^+
\end{array}$                                                                         \\ \hline 
\rowcolor[HTML]{C0C0C0} 
\cellcolor[HTML]{C0C0C0}scenario 7 &  $Ay_{*,0}=\left(\begin{array}{c} c_1^+  \\ -c_3^-+c_5^+ \\ c_3^- \\ -c_3^-+c_1^+ \\ c_5^+\end{array}\right)$ {\color[HTML]{343434} }                         &  $ \def\arraystretch{1.5}
\begin{array}{rcl}
c_2^-\le &-c_3^-+c_5^+&\le c_2^+\\
c_4^-\le &-c_3^-+c_1^+&\le c_4^+
\end{array}$                                          \\ \hline
scenario 8                         & $Ay_{*,0}=\left(\begin{array}{c} c_3^++c_4^+  \\ c_2^+ \\ c_3^+ \\ c_4^+ \\ c_2^++c_3^+\end{array}\right)$ \cellcolor[HTML]{FFFFFF}{\color[HTML]{343434} } & $ \def\arraystretch{1.5}
\begin{array}{rcl}
c_1^-\le &c_3^++c_4^+&\le c_1^+\\
c_5^-\le &c_2^++c_3^+&\le c_5^+
\end{array}$                                            \\ \hline
\end{tabular}
\caption{Same as Table~\ref{table2}, but for scenarios 5-8.}
\label{table3}
\end{table}

The focus of the present section is on the elastoplastic model shown in Fig.~\ref{fig1} (earlier introduced in Rachinskiy \cite{Rachinskiy}), which allows to fully illustrate the practical implementation of Theorem~\ref{mainthm}. 
\begin{figure}[h]
\centering
\includegraphics[scale=0.5]{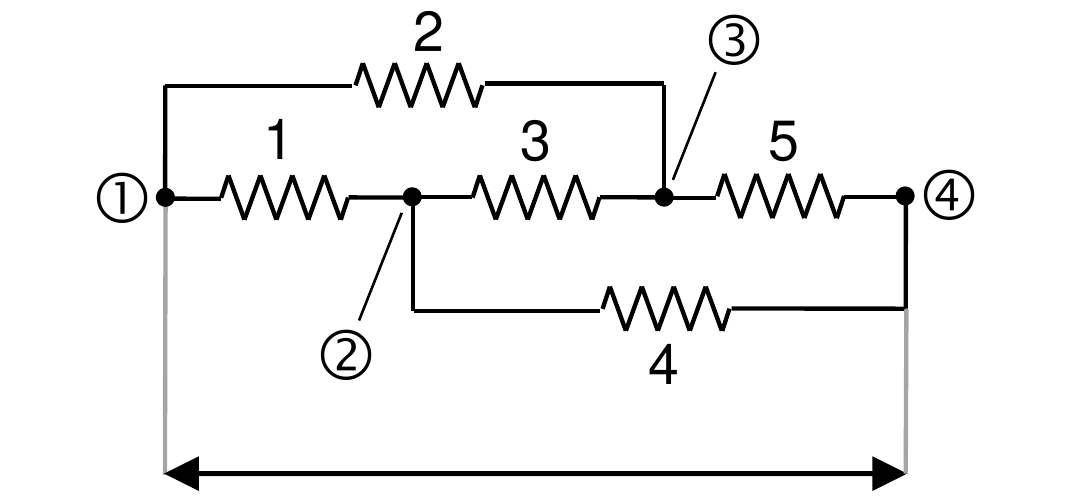}
\caption{A system of 5 elastoplastic springs on 4 nodes that we investigate to illustrate our method. A displacement-controlled loading of gradually increasing or decreasing magnitude is applied as the arrows show. }\label{fig1}
\end{figure}
According to Gudoshnikov-Makarenkov \cite[\S2]{PhysicaD} the elastoplastic system of Fig.~\ref{fig1} leads to the following expressions for  $D$ and $R$
$$
  D=\left(\begin{array}{cccc}
  -1 & 1 & 0 & 0\\
  -1 & 0 & 1 & 0\\
  0 & -1 & 1 & 0\\
  0 & -1 & 0 & 1\\
  0 & 0  & -1 &1
  \end{array}\right),\quad R=\left(\begin{array}{c} 1 \\ 0 \\ 1 \\ 0 \\ 1\end{array}\right).
$$
We now follow Gudoshnikov-Makarenkov \cite[\S5]{PhysicaD} to formulate a sweeping process (\ref{sp}) corresponding to the elastoplastic system $(D,A,C,R,l(t))$. 

\vskip0.2cm

\noindent First of all, based on \cite[formula (17)]{PhysicaD}, we compute the dimension of sweeping process (\ref{sp}) as
$$
   d=m-n+q+1=5-4+1+1=3.
$$
According to \cite[\S5, Step~1]{PhysicaD}, we then look for an $4\times 2$ matrix $M$ such that $R^TDM=0$ and such that the matrix $DM$ is full rank. Such a matrix $M$ can be taken as
$$
  M=\left(\begin{array}{cc}
   0 & 0 \\
   1 & 1 \\
   1 & -1\\
   0 & 0 
   \end{array}\right)\qquad {\rm with}\qquad DM=\left(\begin{array}{cc}
   1 & 1 \\
   1 & -1 \\
   0 & -2\\
   -1 & -1\\
   -1 & 1 
   \end{array}\right).
$$
The next step is determining $\mathcal{V}_{basis}$ which consists of $d=3$ linearly independent columns of $\mathbb{R}^m=\mathbb{R}^5$ and solves $(DM)^TA\mathcal{V}_{basis}=0$. Such a $\mathcal{V}_{basis}$ can be takes as
$$
   \mathcal{V}_{basis}=\left(\begin{array}{ccc} 
   0 & 1/a_1 & 1/a_1\\
   0 & 1/a_2 & -1/a_2\\
   1/a_3 & 0 & 1/a_3\\
   -1/a_4 & 1/a_4 & 0\\
   1/a_5 & 1/a_5 & 0\end{array}\right)\qquad {\rm with}\qquad A\mathcal{V}_{basis}=\left(\begin{array}{ccc}
   0 & 1 & 1\\
   0 & 1 & -1\\
   1 & 0 & 1\\
   -1 & 1 & 0\\
   1 & 1 & 0
   \end{array}\right).
$$

\noindent Finally, a $5\times 2$ full rank matrix $D^\perp$ satisfying $(D^\perp)^TD=0$ can be taken as 
\begin{equation}\label{RTDperpex}
   D^\perp=\left(\begin{array}{cc}
   0 &1 \\
   0 & -1 \\
   1 & 1\\
   -1 & 0\\
   1 & 0
   \end{array}\right)\quad\mbox{leading to}\quad    \left(\begin{array}{c}
  R^T \\
  (D^\perp)^T \end{array}\right)=\left(\begin{array}{ccccc}
  1 & 0 & 1 & 0 & 1\\
  0 & 0 & 1 & -1 & 1\\
  1 & -1 & 1 & 0 & 0
  \end{array}\right).
\end{equation}


\noindent In what follows, we consider 
\begin{equation}\label{dcl}
   l(t)=l_0+l_1\cdot t,\quad \mbox{for some}\ l_0,l_1>0.
\end{equation}

\noindent {\bf Step 1.} We can identify two sets of indexes $I_0$ that we list along with the corresponding inclusion (\ref{step1rev}):
\begin{eqnarray*}
|I_0|=1&:& \ \ \mbox{doesn't work as none of the columns of }\\
&& \ \ \mbox{matrix (\ref{RTDperpex}) are parallel to }(1\  0\  0)^T,\\
I_0=\{(+,1),(+,2)\} &: & \ \  \left(\begin{array}{c} 1 \\ 0 \\ 0 \end{array}\right)\hskip-0.05cm l_1 \hskip0.05cm \in {\rm cone}\left(\left\{\left(\begin{array}{c} -1  \\ 0  \\ -1 \end{array}\right),\left(\begin{array}{c} 0 \\  0 \\  1\end{array}\right)\right\}\right), \\
I_0=\{(+,4),(+,5)\} &: & \ \ \left(\begin{array}{c} 1 \\ 0 \\ 0 \end{array}\right)\hskip-0.05cm l_1 \hskip0.05cm \in{\rm cone}\left(\left\{\left(\begin{array}{c} 0  \\ -1  \\ 0 \end{array}\right),\left(\begin{array}{c} 1 \\  1 \\  0\end{array}\right)\right\}\right), \\
I_0=\{(+,1),(-,3),(+,5)\} &: & \ \ \left(\begin{array}{c} 1 \\ 0 \\ 0 \end{array}\right)\hskip-0.05cm l_1 \hskip0.05cm \in{\rm cone}\left(\left\{\left(\begin{array}{c} 1  \\ 0  \\ 1 \end{array}\right),\left(\begin{array}{c} -1 \\  -1 \\  -1\end{array}\right), \left(\begin{array}{c} 1 \\  1 \\  0\end{array}\right)\right\}\right),   \\
I_0=\{(+,2),(+,3),(+,4)\} &: & \ \ \left(\begin{array}{c} 1 \\ 0 \\ 0 \end{array}\right)\hskip-0.05cm l_1 \hskip0.05cm \in{\rm cone}\left(\left\{\left(\begin{array}{c} 0  \\ 0  \\ -1 \end{array}\right),\left(\begin{array}{c} 1 \\  1 \\  1\end{array}\right), \left(\begin{array}{c} 0 \\  -1 \\  0\end{array}\right)\right\}\right).
\end{eqnarray*}

\noindent {\bf Step 2.} For each $|I_0|<3$ we consider all possible $|I_1|=1$ for which condition (\ref{conditionstep2}) holds. We get a total of three possible $I_1$ for each $I_0$. Each such $I_1$ is extended to $\{I_i\}_{i\in\overline{1,2}}$ according to the procedure described in Section~\ref{guide}. These $\{I_i\}_{i\in\overline{1,2}}$ are listed in column ``Step 2'' of Tables~\ref{table2} and \ref{table3}, thus giving us 8 different scenarios. One or another scenario will take place depending on the feasibility condition that we formulate in the next step.

\vskip0.2cm

\noindent {\bf Step 3.} Fixing $I_0$, $I_1$, and $I_2$ corresponding to scenario~1 of Table~\ref{table1}, we use formula (\ref{y*iformula}) in order to compute $y_{*,1}$ and $y_{*,2}$ as well as to formulate the respective feasibility condition (\ref{y*ifeasible}) which consists of 4 two-sided inequalities. The results of these computations are summarized in line 1 of Table~\ref{table2}. Then we continue by analogy through lines scenarios 2-6 of Table~\ref{table1} and fill out the respective lines of Tables~\ref{table2} and \ref{table3}. For scenario~7 of Table~\ref{table1} formula (\ref{y*iformula}) gives a single vertex $y_{*,0}$ and formula (\ref{y*ifeasible}) gives just a pair of two-sided inequalities that constitute the feasibility condition, see line 3 of Table~\ref{table3}. Computation for the similar scenario 8 are summarized in line 4 of the same table. 

\vskip0.2cm

\noindent Using the fact that each of the inclusions in Step~1 holds in a strict sense (i.e. the vector $(1,0,0)^T$ never belongs to the boundary of the respective cone), we can now use Proposition~\ref{step3prop} to obtain the following  statement about the evolution of the model of Fig.~\ref{fig1}.

\begin{proposition} \label{prop7} Assume that the elastic limits $c^-_i,c_i^+$ of the elastoplastic springs of the model of Fig.~\ref{fig1} with displacement-controlled loading (\ref{dcl}) satisfy the feasibility condition of one of the 8 scenarios of Tables~\ref{table2} and \ref{table3}. 
\begin{itemize}
\item[A:]If the feasibility condition of one of the scenarios 1-6 holds, then there exists an $\eps>0$  such that the 2 springs with the indexes from $I_0$ (of Table~\ref{table1}) undergo plastic deformation for all sufficiently large $t>0$. During this plastic deformation, 
the stresses of the remaining springs can take any constant value from the line segment given by the second column of Tables~\ref{table2} and \ref{table3}.
\item[B:]If the feasibility condition of one of the scenarios 7-8 holds, then there exists an $\eps>0$ such that the 3 springs with the indexes from $I_0$ (of Table~\ref{table1}) undergo plastic deformation for all sufficiently large $t>0$. During this plastic deformation, the stresses of the remaining 2 springs will take the specific constant values given by the second column of Table~\ref{table3}.  
\end{itemize}  
\end{proposition}

\vskip0.1cm

\noindent {\bf Step 4.} {\boldmath{\bf $|I_0|=2.$}} In this case, for any $(\alpha_*,j_*)\in I_0$, the set $I_0\backslash\{(\alpha_*,j_*)\}$ consists of just one element $\{(\alpha,j)\}$ and formula (\ref{newnewnew}) can be rewritten as

\begin{equation}\label{d=2new}
\begin{array}{l}
\eps_0(t)=\min\limits_{(\alpha,j)\in I_0}S_j,\\
S_j=\left\|-c'(t)-{\rm proj}^A(-c'(t),{\rm span}\{n_j\})\right\|^A,\\ {\rm proj}^A(-c'(t),{\rm span}\{n_j\})=-n_j\dfrac{1}{n_j^TAn_j}n_j^T Ac'(t).
\end{array}
\end{equation}

\noindent Therefore, for scenarios 1-6 we get

\noindent $I_0=\{(+,1),(+,2)\}:$ $\eps_0=\min\{S_1,S_2\},$

\noindent $I_0=\{(+,4),(+,5)\}:$ $\eps_0=\min\{S_4,S_5\}.$

\vskip0.2cm

\noindent Computation in Mathematica gives

\vskip0.2cm

$S_1= l_1\sqrt{\dfrac{a_2 \left(a_4 a_5+a_3 \left(a_4+a_5\right)\right)}{a_2 \left(a_3+a_4\right)+a_4 a_5+a_3 \left(a_4+a_5\right)}},$

$
S_2= l_1\sqrt{\dfrac{a_1 \left(a_4 a_5+a_3 \left(a_4+a_5\right)\right)}{a_1 \left(a_3+a_5\right)+a_4 a_5+a_3 \left(a_4+a_5\right)}},$

$
S_4= l_1 \sqrt{\dfrac{a_5\left(a_1 a_2 +a_3(a_1+a_2)\right)}{a_5(a_3+a_1)+a_1 a_2 +a_3(a_1+a_2)}},$

$
S_5= l_1 \sqrt{\dfrac{a_4\left(a_1 a_2 +a_3(a_1+a_2)\right)}{a_4(a_3+a_2)+a_1 a_2 +a_3(a_1+a_2)}},
$

\vskip0.2cm

\noindent which allows to fill out the first two lines of Table~\ref{table4}.

\vskip0.2cm

\noindent {\boldmath{\bf $|I_0|=3.$}} In this case, for any $(\alpha_*,j_*)\in I_0$, the set $I_0\backslash\{(\alpha_*,j_*)\}$ consists of two elements $\{(\alpha_1,j_1),(\alpha_2,j_2)\}$ and formula (\ref{newnewnew}) can be rewritten as
\begin{equation}\label{d>2new} \def\arraystretch{1.5}
\begin{array}{l}
\eps_0(t)=\min\limits_{(\alpha_1,j_1),(\alpha_2,j_2)\in I_0}S_{j_1 j_2},\\
S_{j_1 j_2}=\left\|-c'(t)-{\rm proj}(-c'(t),{\rm span}\{n_{j_1},n_{j_2}\})\right\|^A \\
 {\rm proj}^A(-c'(t),{\rm span}\{n_{j_1},n_{j_2}\})=\\
=-
\left(n_{j_1}\ n_{j_2}\right)\left(\begin{array}{cc} n_{j_1}^TA n_{j_1} & n_{j_1}^TA n_{j_2} \\
n_{j_2}^TAn_{j_1} & n_{j_2}^T An_{j_2}\end{array}\right)^{-1}\left(\begin{array}{c} n_{j_1}^T \\ n_{j_2}^T \end{array}\right)Ac'(t).
\end{array}
\end{equation}

\noindent Therefore, for scenarios 7-8 we get

\noindent $I_0=\{(+,1),(-,3),(+,5)\}:$  $\eps_0=\min\{S_{13},S_{15},S_{35}\}$,
 
\noindent  $I_0=\{(+,2),(+,3),(+,4)\}:$  $\eps_0=\min\{S_{23},S_{24},S_{34}\}$.

\vskip0.2cm

\noindent Computation in Mathematica gives

\vskip0.2cm

$
S_{13}=l_1\sqrt{\dfrac{a_2 a_5}{a_2+a_5}},\ \
S_{35}=l_1\sqrt{\dfrac{a_1 a_4}{a_1+a_4}},\ \
S_{15}=l_1\sqrt{\dfrac{a_2 a_3 a_4}{a_2 a_3+a_2a_4+a_3a_4}},
$

$
S_{23}=l_1\sqrt{\dfrac{a_1 a_4}{a_1+a_4}},\ \
S_{34}=l_1\sqrt{\dfrac{a_2 a_5}{a_2+a_5}},\ \
S_{24}=l_1\sqrt{\dfrac{a_1 a_3 a_5}{a_1 a_3+a_1a_5+a_3a_5}},
$

\vskip0.2cm

\noindent which allows to complete the completion of Table~\ref{table4}.

\vskip0.2cm

\begin{table}[]
\centering
\begin{tabular}{|l|l|}
\hline
                                                                                & Step 4 \\ \hline
\begin{tabular}[c]{@{}l@{}}scenario 1\\ scenario 2\\ scenario 3\end{tabular} &   $
   \eps_0=\min\{h(a_2,a_4),h(a_1,a_5)\}, \   
    h(x,y)=l_1\sqrt{\dfrac{x \left(a_4 a_5+a_3 \left(a_4+a_5\right)\right)}{x \left(a_3+y\right)+a_4 a_5+a_3 \left(a_4+a_5\right)}}   
$     \\ \hline
\begin{tabular}[c]{@{}l@{}}scenario 4\\ scenario 5\\ scenario 6\end{tabular} &  $
\eps_0=\min\{ h(a_5,a_1), h(a_4,a_2)\},\
  h(x,y)=l_1 \sqrt{\dfrac{x\left(a_1 a_2 +a_3(a_1+a_2)\right)}{x(a_3+y)+a_1 a_2 +a_3(a_1+a_2)}}
$
      \\ \hline
scenario 7                                                                     &   $\eps_0=\min\left\{l_1\sqrt{\dfrac{a_2 a_5}{a_2+a_5}},l_1\sqrt{\dfrac{a_1 a_4}{a_1+a_4}},
l_1\sqrt{\dfrac{a_2 a_3 a_4}{a_2 a_3+a_2a_4+a_3a_4}}\right\}$     \\ \hline
scenario 8                                                                     &    $\eps_0=\min\left\{
l_1\sqrt{\dfrac{a_1 a_4}{a_1+a_4}},
l_1\sqrt{\dfrac{a_2 a_5}{a_2+a_5}},
l_1\sqrt{\dfrac{a_1 a_3 a_5}{a_1 a_3+a_1a_5+a_3a_5}}\right\}$    \\ \hline
\end{tabular}
\caption{The distance from the displacement-controlled vector $-c'$ to the boundary of the normal cone formed by the normal vectors with indexes $I_0.$
The notations of the table are introduced and explained in Step~4 of Section~\ref{guide} and computed for model of Fig.~\ref{fig1} in Step~4 of Section~\ref{guide_}.}
\label{table4}
\end{table}

\begin{proposition} {\bf (scenarios 7-8)}\label{prop7-8} Assume that elastic limits of the elastoplastic springs of the model of Fig.~\ref{fig1} with displacement-controlled loading (\ref{dcl}) satisfy the feasibility condition of one the scenarios 7-8 of Table~\ref{table3}. Define $\eps_0$ according to Table~\ref{table4} and put
$$
  \tau_d= \dfrac{1}{\eps_0} \cdot\|A^{-1}c^+-A^{-1}c^-\|^A.
$$
Then, for any initial distribution of stresses, the 3 springs with the indexes from $I_0$ (of Table~\ref{table1}) will undergo plastic deformation for $t\ge\tau_d.$
During this plastic deformation, the stresses of all 5 springs will hold the specific constant values given by the second column of Table~\ref{table3}.  
\end{proposition}


\noindent Proposition~\ref{prop7-8} completes the study of scenarios 7-8 and the next step completes the study of scenarios 1-6.

\vskip0.2cm

\noindent {\boldmath{\bf Step 5. Computing $\sigma_i.$}}  For each of the vertexes $I_1$ and $I_2$ in each of the scenarios 1-6 we setup the matrixes $\mathcal{L}_1$ and $\mathcal{L}_2$ according to formula~(\ref{mathcalL}) and use Mathematica to compute the corrections $\sigma_1$ and $\sigma_2$ as defined by formula~(\ref{eigenvalue}). The results of this computation are summarized in Table~\ref{table5}.

\begin{proposition} {\bf (scenarios 1-6)} \label{prop9} Assume that elastic limits of the elastoplastic springs of the model of Fig.~\ref{fig1} with displacement-controlled loading (\ref{dcl}) satisfy the feasibility condition of one the scenarios 1-6 of Table~\ref{table3}. Define $\eps_0$ according to Table~\ref{table4}, define $\sigma_1,\sigma_2$ according to Table~\ref{table5}, and put
$$
  \tau_d= \dfrac{\max\{\sqrt{\sigma_1},\sqrt{\sigma_2}\}}{\eps_0} \cdot\|A^{-1}c^+-A^{-1}c^-\|^A.
$$
Then, for any initial values of stresses, the 2 springs with the indexes from $I_0$ (of Table~\ref{table1}) will undergo plastic deformation for $t\ge\tau_d.$
During this plastic deformation, 
the stresses of the remaining springs can admit any constant value from the line segment given by the second column of Tables~\ref{table2} and \ref{table3}.
\end{proposition}

\begin{table}[]
\centering
\begin{tabular}{|l|l|}
\hline
            & Step 5                                          \\ \hline
\rowcolor[HTML]{C0C0C0} 
scenario 1 &  $\sigma_1=\sigma_2=\max\left\{1,\dfrac{\left(a_1+a_4\right) \left(a_2+a_5\right) \left(a_4 a_5+a_3 \left(a_4+a_5\right)\right)}{a_4 a_5 \left(a_3 a_4+a_2
\left(a_3+a_4\right)+a_3 a_5+a_4 a_5+a_1 \left(a_2+a_3+a_5\right)\right)}\right\}$                                               \\ \hline
scenario 2 &  $\sigma_1=\sigma_2=\max\left\{1,\dfrac{\left(a_3 \left(a_2+a_5\right)+a_1 \left(a_2+a_3+a_5\right)\right) \left(a_4 a_5+a_3 \left(a_4+a_5\right)\right)}{a_3
a_5 \left(a_3 a_4+a_2 \left(a_3+a_4\right)+a_3 a_5+a_4 a_5+a_1 \left(a_2+a_3+a_5\right)\right)}\right\}$                                               \\ \hline
\rowcolor[HTML]{C0C0C0} 
scenario 3 &  $\sigma_1=\sigma_2=\max\left\{1,\dfrac{\left(a_1 \left(a_2+a_3\right)+a_3 a_4+a_2 \left(a_3+a_4\right)\right) \left(a_4 a_5+a_3 \left(a_4+a_5\right)\right)}{a_3
a_4 \left(a_3 a_4+a_2 \left(a_3+a_4\right)+a_3 a_5+a_4 a_5+a_1 \left(a_2+a_3+a_5\right)\right)}\right\}$                                                \\ \hline
scenario 4 &  $\sigma_1=\sigma_2=\max\left\{1,\dfrac{\left(a_2 a_3+a_1 \left(a_2+a_3\right)\right) \left(a_2 \left(a_3+a_4\right)+a_4 a_5+a_3 \left(a_4+a_5\right)\right)}{a_2
a_3 \left(a_3 a_4+a_2 \left(a_3+a_4\right)+a_3 a_5+a_4 a_5+a_1 \left(a_2+a_3+a_5\right)\right)}\right\}$                                                \\ \hline
\rowcolor[HTML]{C0C0C0} 
scenario 5 &  $\sigma_1=\sigma_2=\max\left\{1,\dfrac{\left(a_2 a_3+a_1 \left(a_2+a_3\right)\right) \left(a_4 a_5+a_1 \left(a_3+a_5\right)+a_3 \left(a_4+a_5\right)\right)}{a_1
a_3 \left(a_3 a_4+a_2 \left(a_3+a_4\right)+a_3 a_5+a_4 a_5+a_1 \left(a_2+a_3+a_5\right)\right)}\right\}$                                                \\ \hline
scenario 6 & $\sigma_1=\sigma_2=\max\left\{1,\dfrac{\left(a_2 a_3+a_1 \left(a_2+a_3\right)\right) \left(a_1+a_4\right) \left(a_2+a_5\right)}{a_1 a_2 \left(a_3 a_4+a_2
\left(a_3+a_4\right)+a_3 a_5+a_4 a_5+a_1 \left(a_2+a_3+a_5\right)\right)}\right\}$                                                                                       \\ \hline
\rowcolor[HTML]{C0C0C0} 
scenario 7 & \cellcolor[HTML]{343434}                        \\ \hline
scenario 8 & \cellcolor[HTML]{343434}{\color[HTML]{656565} } \\ \hline
\end{tabular}
\caption{Corrections for the distance $\eps_0$ computed in Table~\ref{table4}. 
The notations of the table are introduced and explained in Step~5 of Section~\ref{guide} and computed for model of Fig.~\ref{fig1} in Step~5 of Section~\ref{guide_}.}
\label{table5}
\end{table}

\subsection{Remarks}\label{Remarks-Dima}

\noindent {\bf 1.} 
	Under the conditions of Proposition \ref{prop7}, part B, the 3 springs with indexes from $I_0$ undergo the plastic deformation after they saturate
	(see scenarios 7, 8 of Table \ref{table3}). During this plastic deformation,
	the absolute values of the elongations of these springs continue to increase.
	On the other hand, if condition \eqref{A0} is violated, then other terminal states of the same model with increasing displacement-controlled loading (\ref{dcl}) are possible. Table \ref{new} lists such terminal states for the case where $c_i^-=-c_i^+$ for all $i$. In each of these terminal states of stresses, three springs are saturated but only two of them continue to undergo the plastic deformation as $l$ increases. In other words, two springs continue to stretch in the terminal state, while the other three springs maintain a fixed elongation and stress. 

\begin{table}[]
	\centering
	\begin{tabular}{|l|l|l|}
		\hline
		Saturated springs & Springs undergoing plastic & Springs maintaining their\\
		in the terminal state & deformation in the terminal state & length in the terminal state\\ 
		\hline
1, 2, 4 & 1, 2 & 3, 4, 5\\ 
2, 4, 5 & 4, 5 & 1, 2, 3\\ 
1, 2, 3 & 1, 2 & 3, 4, 5\\ 
1, 2, 4 & 1, 2 & 3, 4, 5 \\ 
3, 4, 5 & 4, 5 & 1, 2, 3\\ 
1, 4, 5 & 4, 5 & 1, 2, 3\\ 
		\hline
	\end{tabular}
\caption{Possible terminal states of the model shown in Fig.~\ref{fig1}
with increasing displacement-controlled loading (\ref{dcl}). Each particular terminal state is achieved for a different domain of parameters $c_i^\pm, a_i$ in the parameter state. These scenarios are different from the scenarios listed in Tables \ref{table2}, \ref{table3}. \label{new}}
\end{table}

\medskip
\noindent {\bf 2.} Scenario~7 deserves particular attention. Since the terminal stress of spring~3 is 
	$c_3^-$ (see Table~\ref{table3}), 
	if the initial stress of spring~3 is greater than 
	$c_3^-$, then spring~3 will necessarily contract before it gets to the terminal state of constant stress, even though the entire network of springs stretches (according to formula (\ref{dcl})). That is, the elongation and stress of spring 3 decrease with increasing ``length of the system'' (input) $l$. All the other springs always respond with increasing length to the increasing $l$.
		
		\vskip0.2cm
		
\noindent 	Perhaps even more interestingly,
		there are several scenarios when spring 3 responds non-monotonically to a monotone input $l$. For simplicity, let us assume that $c_i^+=-c_i^-$ for all $i$, hence the maximal absolute value of stress for spring $i$ is 
		$c_i^+$. We will say that a spring {\em saturates} if its stress reaches either the maximal possible value $c_i^+$ or the minimal possible value
		$c_i^-$.
		Let us consider the zero initial state where the elongations of all five springs are zero (all the springs are relaxed) and apply an increasing input $l=l(t)$. The springs can saturate and de-saturate in different order as $l$ increases, depending on the parameters $a_i, c_i^+$. In particular, one can show that the following scenarios with non-monotone behavior of spring 3 are possible.
		If $a_1 a_5>a_2 a_4$, it is easy to see that initially all the springs stretch. Suppose that the first spring to saturate is either spring 1 or spring 5 at a moment $\tau$. Then, after this point, spring 3 will contract.
		On the other hand, if spring 3 is the first one to saturate at a time $\tau_1$, and the second spring to saturate is either spring 1 or spring 5 at a time $\tau_2>\tau_1$, then spring 3 stretches until the moment $\tau_2$ and contracts after this moment. In the latter scenario, spring 3 saturates at the momnet $\tau_1$, undergoes the plastic deformation between the moments $\tau_1$ and $\tau_2$ and de-saturates at the moment $\tau_2$. 
		
		\vskip0.2cm
		
		\noindent Similar examples of a non-monotone response are possible in the complementary case when $a_1 a_5<a_2 a_4$. Here, starting from the zero state,  spring 3 initially contracts, while springs 1, 2, 4, 5 stretch as $l$ increases. If the first to saturate is either spring 2 or spring 4 at a moment $\tau_1$, then spring $3$ stretches for $t>\tau_1$. If spring 3 is the first to staurate at a moment $\tau_1$ and either spring 2 or spring 4 is the second to saturate at a moment $\tau_2>\tau_1$, then spring 3 contracts for $t<\tau_2$ and stretches for $t>\tau_2$, i.e.\ spring 3 de-saturates when another spring (2 or 4) saturates. 

\vskip0.2cm

\noindent {\bf 3.} Some conclusions about finite time stability of the system shown in Fig.~\ref{fig1} can be obtained from the results of
	\cite{Rachinskiy}, which establish the equivalence of 
	a particular class of the sweeping processes and the Prandtl-Ishlinskii model of one-dimensional plasticity, which is one-period stable \cite{PI2,kp,PI1}.

\vskip0.2cm

\noindent	We say that an initial state $y=(y_1,y_2,y_3,y_4,y_5)$ is {\em reachable from zero} if this state can be reached from the zero state $y=0$
	under at least one input $l=l(t)$, $0\le t\le 1$. In the zero state, all  the springs are relaxed, i.e.\ all the stresses are zero. Let us denoty by $\Omega$ the set of all the reachable from zero states.
	One can show that the zero state is reachable from any state $y\in \Omega$.
	Hence any reachable from zero state can be reached from any other reachable from zero state.

\vskip0.2cm

\noindent Let us consider a restriction of the elastoplastic system shown 
	in Fig.~\ref{fig1}, and the corrsponding sweeping process,
	to the set $\Omega$ of reachable from zero states, i.e.\ we consider the initial states from $\Omega$ only.	
Systems shown in Fig.\ \ref{fig1} can be divided into two types depending on their response to an increasing input such as $l(t)=t$, $t\ge 0$, which is applied to the system assuming the zero initial state.  The system will be called {\em anomalous} if the first spring that saturates under such conditions is spring 3 and the second spring that saturates is either spring 1 or spring 5 in the case $a_1 a_5>a_2 a_4$, and either spring 2 or spring 4 in the case $a_1 a_5<a_2 a_4$. According to Remark~2 of this subsection, the response of spring 3 of an anomalous system to increasing inputs is non-monotone. The corresponding solution of the sweeping process makes a transition from one two-dimensional face of the polyhedron to another two-dimensional face at the moment when the second spring saturates and spring 3 simultaneously de-saturates.

\vskip0.2cm

\noindent
The system will be called {\em regular} if it is not anomalous and no two springs saturate simultaneously during the response to an increasing input starting from the zero initial state. Regular and anomalous systems are represented by complementary open domains seprated by their boundary in the parameter space.

\vskip0.2cm

\noindent The next statement follows from the results of Rachinskiy \cite{Rachinskiy}. Its rigorous proof is outside of the scope of the present paper. 
	\begin{proposition}\label{t1}
The restriction of a regular system shown in Fig.\ \ref{fig1} to the class of reachable from zero states
is one period stable under any periodic input $l(t)$.
	\end{proposition}


\section{Conclusions}\label{concl} In this paper we adjusted and applied the ideas of Adly et al \cite{Adly} about finite-time stability of frictional systems to finite-time stability of sweeping processes with polyhedral moving constraints. By using the results obtained we proposed a step-by-step guide to analyze finite-time reachability of plastic deformation in networks of elastoplastic springs. The proposed guide has been tested on a particular example of 5 elastoplastic springs on 4 nodes, and demonstrated that all the required algebraic computations can be executed in Wolfram Mathematica. The Mathematica notebook is uploaded as supplementary material and can be readily used in other networks of elastoplastic springs.

\vskip0.2cm

\noindent Our step-by-step application guide of Section~\ref{guide} addresses a single displacement-controlled loading and a particular way of creating the list of scenarios. However, the finite-time stability results of Section~\ref{sec2} are obtained for general polyhedral sweeping processes (\ref{sp}) and can be applied to arbitrary networks of elastoplastic springs along the lines of Sections~\ref{mechanics} and \ref{guide}. We anticipate that this kind of applications will facilitate collaboration between set-valued analysts and materials scientists.


%
 \section*{Conflict of interest}

 The authors declare that they have no conflict of interest.

\vskip0.5cm

\appendix

\noindent {\large \bf Appendix}

\section{Skipped proofs}\label{skipped}

\noindent {\boldmath{\bf  Proof of implication  (\ref{prop1}) $\Longrightarrow$ (\ref{d-1}).}} By definition (\ref{goodformula2}), if $\xi\in N_C^A(y_{*,i})$, then there exist non-negative numbers $\lambda_1, ..., \lambda_d$ such that
$$
   \xi=(\left\{\alpha n_{j}:(\alpha,j)\in (I_0\cup I_i)\}\right\})(\lambda_1, ..., \lambda_d)^T.
$$ 
But by (\ref{prop1}), $\{\alpha n_{j}:(\alpha,j)\in (I_0\cup I_i)\}$ is a basis of $V$. Therefore, the correspondence between $\xi\in N_C^A(y_{*,i})$ and non-negative 
$\lambda_1, ..., \lambda_d$ is one-to-one. Therefore, any $\xi\in N_C^A(y_{*,i})$ for which the corresponding $\lambda_1, ..., \lambda_d$ contains $\lambda_{i}=0$ is from $\partial N_C^A(y_{*,i})$, which is exactly the statement of formula (\ref{d-1}).\qed

\vskip0.5cm

\noindent {\bf Proof of formula~(\ref{ESAIMestimate}).} By (\ref{C(t)}),
$
\dfrac{1}{a_j}c_j^-\le y_j\le\dfrac{1}{a_j}c_j^+,\ \mbox{for all }y\in C.
$
Therefore,

\noindent $  \max\limits_{u,v\in C}\left(\|u-v\|^A\right)^2=\sum\limits_{j=1}^m a_j(u_j-v_j)^2\le \sum\limits_{j=1}^m \dfrac{1}{a_j}(c_j^+-c_j^-)^2
  = \left(\|A^{-1}c^+-A^{-1}c^-\|^A\right)^2.
$\qed

\vskip0.5cm




\noindent {\bf Proof of the equivalence (\ref{step1rev}) $\Longleftrightarrow$  (\ref{step1conclusion})}. Statement (\ref{step1rev}) implies the existence of $(\lambda_1,...,\lambda_{|I_0|})$ such that
$$
-\mathcal{V}_{basis}W^{-1}\left(\begin{array}{c}
1\\
0_{m-n+1}\end{array}\right)l'(t)= \mathcal{V}_{basis}W^{-1}\left\{
\left(\begin{array}{c}
  R^T \\
  (D^\perp)^T \end{array}\right)\left(\left\{\alpha e_j:(\alpha,j)\in I_0\right\}\right)\right\}\left(\begin{array}{c}
  \lambda_1\\ \vdots \\ \lambda_{|I_0|}\end{array}\right). 
$$
Due to (\ref{barL}), (\ref{c(t)0}), and (\ref{nj}), the latter formula just coincides with
$$
   c'(t)=\left\{\alpha n_j: (\alpha,j)\in I_0\right\}\left(\begin{array}{c}
  \lambda_1\\ \vdots \\ \lambda_{|I_0|}\end{array}\right),
$$
which yields (\ref{step1conclusion}).\qed

\vskip0.5cm

\noindent {\bf Proof of the implication (\ref{F(t)})-(\ref{Ffeasible}) $\Longrightarrow$ (\ref{y*iformula}).} Based on formula (\ref{prop1}), finding $y_{*,i}$ means solving a system of $d$ algebraic equations
$$
   \left<e_j,A y_{*,i}\right>=c_j^\alpha,\quad (\alpha,j)\in I_0\cup I_i,
$$
or, equivalently,
$$
 \left(\left\{e_j,(\alpha,j)\in I_0\cup I_i\right\}\right)^T  A \mathcal{V}_{basis} v_{*,i}= \left(\left\{c_j^\alpha,(\alpha,j)\in I_0\cup I_i\right\}\right)^T,
$$
where
$
   y_{*,i}=\mathcal{V}_{basis} v_{*,i}.
$\qed

\vskip0.5cm

\section{Technical lemmas} \label{appB}

\begin{lemma}\label{techlemma0} If a non-negative continuously differentiable function $v(t)$ satisfies the differential inequality $\dot v(t)\le -2\eps\sqrt{v(t)}$, then $v(t_1)=0$ for some $t_1\le\dfrac{1}{\eps}v(0).$
\end{lemma}

\noindent {\bf Proof.} The proof follows by observing that the solution of the differential equation $\dot {\bar v}(t)=-2\eps\sqrt{\bar v(t)}$ with $\bar v(0)\ge 0$ is given by
$\bar v(t)=\left(-\eps t+\sqrt{\bar v(0)}\right)^2$ on $[0,\bar t_1],$ where $\bar t_1=\dfrac{1}{\eps}\sqrt{\bar v(0)}.$ \qed

\begin{lemma} \label{techlemma1} Consider $f,g: \mathcal{V}\to \mathcal{V}_1,$ where $\mathcal{V},$ $\mathcal{V}_1$ are scalar product spaces.  If both $D_\xi f(v)$ and $D_\xi g(v)$ exist then $D_\xi\left<f(\cdot),g(\cdot)\right>(v)$ exists and
$$
  D_\xi\left<f(\cdot),g(\cdot)\right>(v)=\left<D_\xi f(v),g(v)\right>+\left<f(v),D_\xi g(v)\right>.
$$
\end{lemma}

\noindent {\bf Proof.} We have
\begin{eqnarray*}
D_\xi\left<f(\cdot),g(\cdot)\right>(v)&=&\lim_{\tau\to 0}\dfrac{\left<f(v+\tau\xi),g(v+\tau\xi)\right>-\left<f(v),g(v)\right>}{\tau}=\\
&=& \left<\lim_{\tau\to 0}\dfrac{f(v+\tau\xi)-f(v)}{\tau},g(v)\right>+\left<f(v),\lim_{\tau\to 0}\dfrac{g(v+\tau\xi)-g(v)}{\tau}\right>+\\
&&+\lim_{\tau\to 0}\left<f(v+\tau\xi)-f(v),\dfrac{g(v+\tau\xi)-g(v)}{\tau}\right>=\\
&=& \left<D_\xi f(v),g(v)\right>+\left<f(v),D_\xi g(v)\right>,
\end{eqnarray*}
where we used that
$$
  \left|\left<f(v+\tau\xi)-f(v),\dfrac{g(v+\tau\xi)-g(v)}{\tau}\right>\right|^2\le\tau\cdot\left\|\dfrac{f(v+\tau\xi)-f(v)}{\tau}\right\|\cdot\left\|\dfrac{g(v+\tau\xi)-g(v)}{\tau}\right\|
$$
by Cauchy-Schwartz inequality.\qed

\begin{lemma} \label{techlemma2} Consider $f: \mathcal{V}\to V_1$ and $u:\mathbb{R}\to \mathcal{V}$, where $\mathcal{V}$, $\mathcal{V}_1$ are scalar product spaces. If both $u'(t_0)$ and $(f\circ u)'(t_0)$ exist and if $f$ is Lipschitz continuous in the neighborhood of $u_0=u(t_0)$, then $D_{u'(t_0)}f(u_0)$ exists and
$$
    D_{u'(t_0)}f(u_0)=(f\circ u)'(t_0).
$$
\end{lemma}

\noindent {\bf Proof.} We have
\begin{eqnarray*}
 && D_{u'(t_0)}f(u_0)=\lim_{\tau\to 0}\dfrac{f(u_0+\tau u'(t_0))-f(u_0)}{\tau}=\\
 && \quad = \lim_{\tau\to 0}\left(\dfrac{f(u(t_0)+\tau u'(t_0))-f(u(t_0+\tau))}{\tau}+\dfrac{f(u(t_0+\tau))-f(u_0)}{\tau}\right)=(f\circ u)'(t_0),
\end{eqnarray*}
where we used Lipschitz continuity of $f$ to conclude that the first fraction in the limit converges to 0 as $\tau\to 0.$\qed

\begin{lemma}\label{centrallemma} If 
$D_\xi{\rm proj}^A(\cdot,F)(v)$ exists then formula (\ref{V1}) holds. 
\end{lemma}

\noindent {\bf Proof.} Let $v,\xi\in\mathbb{R}^d$ be such that 
$D_\xi{\rm proj}^A(\cdot,F)(v)$ exists.
We claim that
\begin{equation}\label{weclaim}
\left<D_\xi{\rm proj}^A(\cdot,F)(v),A(v-{\rm proj}^A(v,F))\right>=0.
\end{equation}
Assume that (\ref{weclaim}) doesn't hold.

\vskip0.2cm

\noindent {\bf Case 1:} Assume that $\left<D_\xi{\rm proj}^A(\cdot,F)(v),A(v-{\rm proj}^A(v,F))\right>>0.$
By the definition of the bilateral directional derivative, 
$$
D_\xi{\rm proj}^A(\cdot,F)(v)=   \dfrac{{\rm proj}^A(v+\tau\xi,F)-{\rm proj}^A(v,F)}{\tau}+\dfrac{o(\tau)}{\tau},
$$
for all $\tau\in\mathbb{R}$ with $|\tau|$ sufficiently small. Therefore, we can choose a sufficiently small positive $\tau$ such that 
$$\left<\dfrac{{\rm proj}^A(v+\tau\xi,F)-{\rm proj}^A(v,F)}{\tau},A(v-{\rm proj}^A(v,F))\right>>0,$$
or, by multiplying by $\tau$,
\begin{equation}\label{formula3}
\left<{\rm proj}^A(v+\tau\xi,F)-{\rm proj}^A(v,F),A(v-{\rm proj}^A(v,F))\right>>0.
\end{equation}
This contradicts the following property of the projection ${\rm proj}^A(v,F)$ (see e.g. Bauschke-Combettes \cite[Theorem~3.16]{Bauschke-Combettes}):
$$
   \left<v_1-{\rm proj}^A(v,F),A(v-{\rm proj}^A(v,F))\right>\le 0, \quad {\rm for\ all\ }v_1\in F.
$$

\noindent {\bf Case 2:} Assume that $\left<D_\xi{\rm proj}^A(\cdot,F)(v),A(v-{\rm proj}^A(v,F))\right><0.$ In this case, we will choose a negative $\tau$ with sufficiently small absolute value $|\tau|$ so that 
$$\left<\dfrac{{\rm proj}^A(v+\tau\xi,F)-{\rm proj}^A(v,F)}{\tau},A(v-{\rm proj}^A(v,F))\right><0,$$ which leads to the same (\ref{formula3}) upon multiplying by $\tau$. 

\vskip0.2cm

\noindent The proof of the lemma is complete.\qed

\begin{lemma}\label{W} For $m \ge n$, consider a $m\times n$-matrix $D$ and $m\times(m-n+1)$-matrix $D^\perp$, such that $(D^\perp)^TD=0_{n\times(m-n+1)}$. If (\ref{rankD}) and (\ref{rankDperp}) hold, then 
\begin{equation}\label{temp0}
   D^\perp\mathbb{R}^{m-n+1}=(D\mathbb{R}^n)^\perp.
\end{equation}
\end{lemma}

\noindent {\bf Proof.}  By the definition of $D^\perp$, 
\begin{equation}\label{temp1}
D^\perp \mathbb{R}^{m-n+1}\subset  {\rm Ker} D^T.
\end{equation}  
Furthermore, we have  
\begin{equation}\label{temp2}
(D\mathbb{R}^n)^\perp={\rm Ker}\hskip0.05cm D^T,
\end{equation} see e.g. Friedberg et al \cite[Exercise 17, p. 367]{Friedberg}. To prove the backwards implication in (\ref{temp0}), we use (\ref{temp2}) and assumption 
(\ref{rankD}) to conclude that 
${\rm Ker}\hskip0.05cm D^T={\rm dim}\left((D\mathbb{R}^n)^\perp\right)=m-n+1$. On the other hand, assumption (\ref{rankDperp})  implies that ${\rm dim}\left(D^\perp\mathbb{R}^{m-n+1}\right)=m-n+1$ too. Therefore, the dimensions of the spaces in the two sides of (\ref{temp1}) coincide and the inclusion (\ref{temp1}) is actually an equality. \qed

\begin{corollary}\label{cor2} Assume that $m\ge n.$ Let $R$ be an $m\times q$-matrix. Let $D^\perp$ be as defined in Lemma~\ref{W}.
Consider
$$
\mathcal{U}=\left\{x\in D\mathbb{R}^n:R^T x=0\right\}.
$$
If conditions (\ref{rankD}) and (\ref{rankDperp}) hold, then
$$  x\in \mathcal{U}\quad\mbox{if and only if}\quad \left(\begin{array}{c}
  R^T\\
  (D^\perp)^T\end{array}\right)x=0.
  $$
\end{corollary}

\noindent {\bf Proof.} The proof follows by observing that $(D^\perp)^Tx=0$ if and only if 
$$
  x\in{\rm Ker}\left((D^\perp)^T\right)=\left(D^\perp\mathbb{R}^{m-n+1}\right)^\perp=\left((D\mathbb{R}^n)^\perp\right)^\perp=D\mathbb{R}^n,
$$
where the first equality is the property that we already used in the proof of Lemma~\ref{W} (see formula (\ref{temp2})) and the second equality is the conclusion of Lemma~\ref{W}.\qed

\begin{corollary}\label{corrank} In the settings of Corollary~\ref{cor2}, assume that 
\begin{equation}\label{rankDTR}
{\rm rank}(D^T R)=q,
\end{equation}
in addition to (\ref{rankD}) and (\ref{rankDperp}). Put $d=m-n+q+1$. Let $\mathcal{V}_{basis}$ be a matrix of $d$ linearly independent vectors of $\mathbb{R}^m$ which are orthogonal to vectors of $\mathcal{U}$ in some scalar product. Then, 
\begin{itemize}
\item[(i)] the $d\times d$-matrix 
$
   \left(\begin{array}{c} R^T \\ (D^\perp)^T\end{array}\right)\mathcal{V}_{basis}
$
is invertible,
\item[(ii)] ${\rm rank}
   \left(\begin{array}{c} R^T \\ (D^\perp)^T\end{array}\right)=m-n+q+1.
$
\end{itemize} 
\end{corollary}

\noindent {\bf Proof.} (i) If 
   $\left(\begin{array}{c} R^T \\ (D^\perp)^T\end{array}\right)\mathcal{V}_{basis}v=0$ for some $v\in\mathbb{R}^d$, then $\mathcal{V}_{basis}v$ must be an element of $\mathcal{U}$ by 
 Corollary~\ref{cor2}. On the other hand,   
  vector $\mathcal{V}_{basis}v$ is orthogonal to the vectors of $\mathcal{U}$, which implies $\mathcal{V}_{basis}v=0$ which can only happen if $v=0.$
  
\vskip0.2cm

\noindent (ii) By the rank-nullity theorem (see e.g. Friedberg et al  \cite[Theorem~2.3]{Friedberg}) and by Corollary~\ref{cor2} we have

\vskip0.2cm

$
{\rm rank}
   \left(\begin{array}{c} R^T \\ (D^\perp)^T\end{array}\right)=
   m-{\rm dim}\left({\rm ker}\left(\begin{array}{c} R^T \\ (D^\perp)^T\end{array}\right)\right)=m-{\rm dim}(\mathcal{U}).
 $ 

\vskip0.2cm

\noindent In this formula, ${\rm dim}(\mathcal{U})=n-q-1$ by Gudoshnikov-Makarenkov \cite[Lemma~3.8]{ESAIM}.
  \qed

\begin{lemma} \mbox{\rm (Rockafellar-Wets \cite[Theorem~6.46]{Rockafellar-Wets})} \label{Rockafellar-Wets-Theorem6.46}
Consider a polyhedron 
$$
   C=\bigcap_{k=1}^K\left\{v\in\mathbb{R}^d:\left<n_k,v\right>\le c_k\right\},
$$
where $n_k\in\mathbb{R}^d$, $c_k\in\mathbb{R},$ $K\in\mathbb{N}.$ If $I(v)=\left\{k\in\overline{1,K}:\left<n_k,v\right>=c_k\right\}$, then
$$
   N_\mathcal{C}(y)={\rm cone}\left\{n_k:k\in I(v)\right\}.
$$
\end{lemma}

\noindent {\bf Proof of Lemma~\ref{N_C}.} Fix $y\in \mathcal{V}.$ The definition of $N_{\widetilde C}^A(y)$ reads as 
\begin{equation}\label{readsas1}
  \left<N_{\widetilde C}^A(y),A(\widetilde c-y)\right>\le 0,\quad \widetilde c\in {\widetilde C}.
\end{equation}
Let $d$ be the dimension of $\mathcal{V}$ and let $\mathcal{V}_{basis}$ be a $m\times d$-matrix of some linearly independent vectors of $\mathcal{V}$. Then we can represent $\widetilde C$ as
$$
\widetilde C=\mathcal{V}_{basis}C, \quad {\rm where}\ C=\bigcap_{k=1}^K\left\{v\in\mathbb{R}^d:\left<n_k,v\right>\le c_k\right\},\ n_k=(A\mathcal{V}_{basis})^T\widetilde n_k.
$$
Defining $v\in\mathbb{R}^d$ in such a way that $y=\mathcal{V}_{basis}v$, statement (\ref{readsas1}) can be rewritten as
$$
  \left<N_{\widetilde C}^A(\mathcal{V}_{basis}v),A(\widetilde c-\mathcal{V}_{basis}v)\right>\le 0,\quad \widetilde c\in \mathcal{V}_{basis}C,
$$
or
$$
  \left<(A\mathcal{V}_{basis})^TN_{\widetilde C}^A(\mathcal{V}_{basis}v),c-v\right>\le 0,\quad  c\in C.
$$
But the definition of $N_C(v)$ reads as
$$
  \left<N_C(v),c-v\right>\le 0,\quad c\in C.
$$
Therefore, $(A\mathcal{V}_{basis})^TN_{\widetilde C}^A(\mathcal{V}_{basis}v)=N_C(v)$ or, incorporating the conclusion of Lemma~\ref{Rockafellar-Wets-Theorem6.46},
$$
  (A\mathcal{V}_{basis})^TN_{\widetilde C}^A(\mathcal{V}_{basis}v)={\rm cone}\left\{(A\mathcal{V}_{basis})^T\widetilde n_k:k\in I(v)\right\},
$$
from where the required statement follows.\qed

\begin{proposition}\label{proj} For any convex set $F\subset\mathbb{R}^m$, 
$${\rm proj}^A(v,F)+c={\rm proj}^A(v+c,F+c),\quad v,c\in F.$$
\end{proposition}

\noindent {\bf Proof.} Indeed, let
$$
  v_*''={\rm proj}^A(v+c,F+c).
$$
Then $v_*''$ satisfies
$$
   \min_{v''\in F+c}\|v+c-v''\|^A=\|v+c-v_*''\|^A,
$$
or 
$$
   \|v+c-v''\|^A>\|v+c-v_*''\|^A,\quad\mbox{for all }v''\in F+c,\ v''\not=v_*'',
$$
or
$$
   \|v-v'\|^A>\|v+c-v_*''\|^A,\quad\mbox{for all }v'\in F,\ c+v'\not=v_*''.
$$
Introducing $v_*'=v_*''-c$, 
$$
   \|v-v'\|^A>\|v-v_*'\|^A,\quad\mbox{for all }v'\in F,\ v'\not=v_*'.
$$
Therefore,
$$
   \min_{v'\in F}\|v-v'\|^A=\|v-v_*'\|^A,
$$
i.e.
$
   v_*'={\rm proj}(v,F).
$\qed

\vskip0.5cm

\noindent {\bf Proof of Lemma~\ref{computedistance}.} Invertibility of the $k\times k$-matrix $\mathcal{N}^T A \mathcal{N}$ follows from the fact that ${\rm rank}(\sqrt{A}\mathcal{N})=k$ and so ${\rm rank}(\mathcal{N}^T A \mathcal{N})={\rm rank}((\sqrt{A}\mathcal{N})^T\sqrt{A}\mathcal{N})=k$, see e.g. Friedberg et al \cite[\S6.3, Lemma~2]{Friedberg}.

\vskip0.2cm

\noindent To prove formula (\ref{distcone}), we observe that 
$$
   {\rm dist}^A\left(-c', {\rm span}\left\{n_{i_1},...,n_{i_k}\right\}\right)=\left\|-c'-{\rm proj}^A\left(-c',{\rm span}\left\{n_{i_1},...,n_{i_k}\right\}\right)\right\|^A.
$$
By the definition of projection (see e.g. Bauschke-Combettes \cite[\S3.2]{Bauschke-Combettes}), 
$${\rm proj}^A\left(-c',{\rm span}\left\{n_{i_1},...,n_{i_k}\right\}\right)=\lambda_1 n_{i_1}+\ldots+\lambda_k n_{i_k},$$ where $\lambda_1,\ldots,\lambda_k \in\mathbb{R}$ minimize the quantity
$$
   \left<-c'-\lambda_1 n_{i_1}-\ldots-\lambda_k,A(-c'-\lambda_1 n_{i_1}-\ldots-\lambda_k)\right>.
$$ 
Therefore, 
$$
\begin{array}{rcl}
 \left<-c'-\lambda_1 n_{i_1}-\ldots-\lambda_k,An_{i_1}\right>&=&0,\\
& \vdots&\\
 \left<-c'-\lambda_1 n_{i_1}-\ldots-\lambda_k,An_{i_k}\right>&=&0,
\end{array}
$$
for the unknown $\lambda_1,...,\lambda_k,$ or, equivalently,
$$
  -\mathcal{N}^T A c'-\mathcal{N}^T A \mathcal{N}(\lambda_1\ldots \lambda_k)^T=0.
$$
Formula (\ref{distcone}) follows by solving this equation for $(\lambda_1\ldots \lambda_k)^T$ and by plugging the result into
$
  {\rm proj}^A\left(-c',{\rm span}\left\{n_{i_1},...,n_{i_k}\right\}\right)=\mathcal{N}(\lambda_1\ldots \lambda_k)^T.$ 
 \qed

\vskip0.2cm

\begin{lemma}\label{enumeration} If conditions (\ref{prop1}), (\ref{noother}), and (\ref{independent}) hold, then all vertices of $F$ are contained in the set $\{y_{*,1},...,y_{*,M}\}.$
\end{lemma}

\noindent {\bf Proof.} Assume that $F$ has a vertex $\tilde y_*\not\in \{y_{*,1},...,y_{*,M}\}.$ We have
$$
   \{\tilde y_*\}=\{y:y\in\overline{\overline{L}}(\alpha,j),\ (\alpha,j)\in I_0\cup \{j_1,...,j_{d-|I_0|}\}\},
$$
where $|I_0\cup\{j_1,...,j_{d-|I_0|}\}|=d.$ By (\ref{noother}), 
$$
   \{j_1,...,j_{d-|I_0|}\}=\cup\bigcap_{i\in J_{\tilde y_*}} I_i.
$$
But $|I_i|=d-|I_0|$ by (\ref{independent}). Therefore, there exists $\tilde i\in J_{\tilde y_*}$ such that $\{j_1,...,j_{d-|I_0|}\}=I_{\tilde i},$ i.e. $\tilde y_*=y_{*,\tilde i}.$ The proof of the lemma is complete.\qed



\end{document}